\documentclass[11pt]{amsart}
\usepackage{amsmath, amsmath,amssymb,latexsym,dsfont}
\usepackage[left=2.6cm,right=2.6cm,top=2.7cm,bottom=2.7cm]{geometry}
\usepackage{amsfonts,epsfig, textcomp,esint}
\usepackage{graphicx}
\usepackage{cases}
\usepackage{color}
\usepackage{tikz}
\usepackage{cancel}
\usepackage{amscd}
\usetikzlibrary{positioning}
\usepackage[french,english]{babel}
\usepackage{subfig}
\usepackage[OT1]{fontenc}
\usepackage{amsthm}
\usepackage{enumerate}
\usepackage{hyperref}
\usepackage{cleveref}

\theoremstyle{plain}
\newtheorem{theorem}{Theorem}[section]
\newtheorem{lemma}{Lemma}[section]

\newtheorem{proposition}{Proposition}[section]
\newtheorem{corollary}{Corollary}[section]
\newtheorem{remark}{Remark}[section]

\newtheorem{definition}{Definition}[section]
\numberwithin{equation}{section}

\newcommand{\tmu}{\widetilde \mu}

\newcommand{\hy}{\hat y}
\newcommand{\hu}{\hat u}

\newcommand{\hw}{\hat w}

\newcommand{\cH}{\mathcal H}

\newcommand{\cG}{\mathcal G}

\newcommand{\tlambda}{\widetilde \lambda}

\newcommand{\ty}{\widetilde y}

\newcommand{\cN}{{\mathcal N}}

\newcommand{\N}{\mathbb{N}}

\newcommand{\loc}{\operatorname{loc}}

\newcommand{\eps}{\varepsilon}
\newcommand{\mR}{\mathbb{R}}
\newcommand{\mZ}{\mathbb{Z}}

\newcommand{\mC}{\mathbb{C}}

\newcommand{\supp}{\operatorname{supp}}
\newcommand{\dsp}{\displaystyle}

\newcommand{\M}{{\mathcal M}}

\newcommand{\diff}[1][-3]{\mathop{}\mkern#1mu{d}}
\numberwithin{equation}{section}

\newcounter{proofstep}

\let\oldproof\proof
\renewcommand{\proof}{\oldproof\setcounter{proofstep}{0}}

      \makeatletter
      \def\@setcopyright{}
      \def\serieslogo@{}
      \makeatother

\usepackage[normalem]{ ulem }
\usepackage{soul}

\definecolor{darkgreen}{rgb}{0,0.7,0}

 \newcommand{\be}{\begin{equation}}
\newcommand{\ee}{\end{equation}}

\newcommand{\dy}{\delta y}
\newcommand{\cM}{\mathcal M}

\title[Local controllability in finite time and the controllable time of KdV]{Local controllability in finite time and the controllable time  of the Korteweg-De Vries equation using the right Neumann controls}

 \author{Hoai-Minh Nguyen}
\address[Hoai-Minh Nguyen]{Sorbonne Universit\'e, Universits\'e Paris Cit\'e, CNRS, INRIA, Laboratoire Jacques-Louis Lions, LJLL, F-75005 Paris, France}
\email{hoai-minh.nguyen@sorbonne-universite.fr}

\begin{document}

\begin{abstract} We investigate the local boundary controllability of the Korteweg–de Vries (KdV) equation with right Neumann boundary controls at critical lengths. We show that the KdV system is not locally null-controllable in small time for all critical lengths for which the unreachable subspace of the linearized system has dimension at least two. This result extends the work of Coron, Koenig, and Nguyen, who established it for a subclass  of these lengths. We also obtain a new controllability time for such systems for all but two critical lengths. It is worth noting that the latest results on the controllability time prior to this work date back to the work of Cerpa (2007) and  Cerpa and Crépéau (2009).  
\end{abstract}

\maketitle



\tableofcontents

\section{Introduction}

This paper is devoted to the local controllability properties of the Korteweg-de Vries (KdV) equation using the right Neuman control. More precisely, we consider the following control problem, for $T>0$,  
\begin{equation}\label{sys-KdV-N}\left\{
\begin{array}{cl}
y_t (t, x) + y_x (t, x)  + y_{xxx} (t, x)  + y (t, x)  y_x (t, x)  = 0 &  \mbox{ in } (0, T) \times (0, L), \\[6pt]
y(\cdot, 0) = y(\cdot, L) = 0 & \mbox{ in }   (0, T), \\[6pt]
y_x(\cdot, L) = u & \mbox{ in }  (0, T), \\[6pt]
y(0, \cdot)  = y_0  &  \mbox{ in }  (0, L). 
\end{array}\right.
\end{equation}
Here  $y$ is the state, $y_0 \in L^2(0, L)$ is an initial datum,  and $u$ is a control, belonging to an appropriate functional space.  

The KdV equation has been introduced by Boussinesq \cite{1877-Boussinesq} and Korteweg
and de Vries \cite{KdV} as a model for the propagation of surface water waves along a
channel. This equation also furnishes a very useful nonlinear approximation model
including a balance between weak nonlinearity and weak dispersive effects, see, e.g., \cite{Whitham74, Miura76, Kato83}.  The KdV
equation has been investigated from various aspects of mathematics, including the well-posedness, the existence, and stability of solitary waves, the integrability, the
long-time behavior,  see,  e.g., ~\cite{Whitham74, Miura76, Kato83, Tao06, LP15}.

The local controllability for the KdV equation has been studied extensively in the literature, see, e.g., the surveys \cite{RZ09, Cerpa14} and the references therein. We briefly review here some results concerning the boundary controls.  When the controls are $y(\cdot, 0)$, $y(\cdot, L)$,   $y_x(\cdot, L)$, Russell and Zhang \cite{RZ96} proved that the KdV equation is small-time,  locally exactly controllable. When the controls are  $y(\cdot, 0)$ and $y(\cdot, L)$ ($y_x(\cdot, L) = 0$), Glass and Guerrero \cite{GG08} showed that the system is locally exactly controllable. 
The case of the left boundary control ($y(\cdot, L) = y_x(\cdot, L) = 0$) was investigated by Rosier \cite{Rosier04} (see also \cite{GG08}). The small-time local, null controllability holds in this case. The exact controllability does not hold for initial and final datum in $L^2(0, L)$ due to the regularization effect, but 
holds for a subclass of infinitely smooth initial and final data \cite{MRRR19}.

The controllability properties of system \eqref{sys-KdV-N} have been studied
intensively. For initial and final datum in $L^2(0, L)$ and controls in $L^2(0, T)$,
Rosier~\cite{Rosier97} proved that the system
is small-time locally controllable around 0 provided that the length $L$ is not critical,
i.e., $L \notin \cN_N$, where
\begin{equation}\label{def-cN}
\cN_N : = \left\{ 2 \pi \sqrt{\frac{k^2 + kl + l^2}{3}}; \, k, l \in \N_*\right\}.
\end{equation}
To this end, he studied the controllability of the linearized system using the Hilbert
uniqueness method  and compactness-uniqueness arguments, and  showed that the
linearized system is controllable if $L \not \in \cN_N$. He as well established
that when $L \in \cN_N$,  the linearized system is not controllable. More precisely, Rosier showed
that there exists a non-trivial, finite-dimensional subspace $\M_N$ of $L^2(0, L)$  such that its orthogonal space is reachable from $0$ whereas  $\M_N$ is not (see, e.g.,  \eqref{def-MN} for the definition of $\cM_N$).

To tackle the control problem for the critical length $L \in \cN_N$ with  initial and final  datum in $L^2(0, L)$ and controls in $L^2(0, T)$, Coron and Cr\'epeau  \cite{CC04}
introduced the
power series expansion method. The idea is to take into account the effect of
 the nonlinear term $y y_x$  absent in  the linearized system. Using this
method, they showed \cite{CC04} (see also \cite[section 8.2]{Coron07}) that  system \eqref{sys-KdV-N} is
small-time locally controllable if $L = m 2 \pi$ for $m \in \N_*$ satisfying
\begin{equation}
\nexists (k, l) \in \N_* \times \N_* \mbox{ with } k^2 + kl + l^2 = 3 m^2 \mbox{ and }
k \neq l.
\end{equation}
In this case,  $\dim \M_N = 1$ and $\M_N$ is spanned by $1 - \cos x$. Cerpa \cite{Cerpa07}
developed the analysis in \cite{CC04} to prove that system \eqref{sys-KdV-N}  is locally
controllable in \emph{finite time} in the case $\dim \M_N = 2$. This corresponds to
the
case where
\be \label{def-L}
L = 2 \pi \sqrt{\frac{k^2 + kl + l^2}{3}},
\ee
for some $k, \,  l \in \N_*$ with   $k>l$, and there is no $m, n \in N_*$ with $m>n$
and $m^2 + mn + n^2 = k^2 + kl + l^2$. Later, Cr\'epeau and Cerpa \cite{CC09}
succeeded to extend the ideas in \cite{Cerpa07} to obtain the local
controllability for all other critical lengths in {\it finite time} as well. To summarize,
concerning the critical lengths
with  initial and final  datum in $L^2(0, L)$ and controls in $L^2(0, T)$, the {\it small-time} local controllability is valid when  $\dim \M_N
= 1$ and the local controllability in {\it finite time} holds when $\dim \M_N \ge 2$.

\subsection{Statement of the main results and ideas of the proof}

With Coron and Koenig \cite{CKN-24}, we recently studied 
the small-time controllability property of system \eqref{sys-KdV-N} when $\dim \cM_N \ge 2$. In this direction,  we proved that the control system \eqref{sys-KdV-N} is {\it not} small-time controllable using controls $u$ in $H^1(0, T)$ if \eqref{def-L} for some $k, l \in \N_*$ such that $2k + l \not \in 3 N_*$. 

Our approach in \cite{CKN-24} is inspired by the power
series expansion method introduced by Coron and Cr\'epeau  \cite{CC04}. The idea of this method is to search/understand a control $u$ of the form
\[
u = \eps u_1 + \eps^2 u_2 +  \cdots.
\]
The corresponding solution then formally takes the form
\[
y = \eps y_1 + \eps^2 y_2  + \cdots,
\]
and the non-linear term $y y_x$ can be written as
\[
y y_x = \eps^2 y_1 y_{1, x} + \cdots.
\]
One then obtains the following systems for $y_1$ and $y_2$: 
\begin{equation}\label{eq:first_order}\left\{
\begin{array}{cl}
y_{1, t} (t, x) + y_{1, x} (t, x) + y_{1, xxx} (t, x) = 0 &  \mbox{ for } t \in
(0, T), \, x \in (0, L), \\[6pt]
y_1(t, x=0) = y_1(t, x=L) = 0 & \mbox{ for } t \in (0, T), \\[6pt]
y_{1, x}(t , x= L) = u_1(t) & \mbox{ for } t \in (0, T),
\end{array}\right.
\end{equation}
\begin{equation} \label{eq:second_order}\left\{
\begin{array}{cl}
y_{2, t} (t, x) + y_{2, x} (t, x) + y_{2, xxx} (t, x) + y_1 (t, x) y_{1, x} (t, x)  = 0 &
\mbox{ for } t \in (0, T), \, x \in (0, L), \\[6pt]
y_2(t, x=0) = y_2(t, x=L)  = 0 & \mbox{ for } t \in (0, T), \\[6pt]
y_{2, x}(t , x= L) = u_2(t) & \mbox{ for } t \in (0, T). 
\end{array}\right.
\end{equation}

Assume that \eqref{def-L} holds for some $k, l \in \N_*$. Set 
\begin{equation}\label{def-p}
p = \frac{(2k + l)(k-l)(2 l + k)}{3 \sqrt{3}(k^2 + kl + l^2)^{3/2}}, 
\end{equation}
and
\begin{equation}\label{eta-kdv}
\eta_1 = - \frac{2 \pi i}{3 L} (2k + l), \quad \eta_2 = \eta_1 + \frac{2\pi i}{L} k, \quad
\eta_3 = \eta_2 + \frac{2\pi i}{L} l. 
\end{equation}
Then $\eta_j$ with $1 \le j \le 3$ are the solutions of the equation 
\be \label{eta-p}
\eta^3 + \eta - i p =0. 
\ee

Denote \footnote{Here and in what follows, we use the convenction $\eta_{j+3} = \eta_j$ for $j \ge 1$.}
\begin{equation}\label{def-psi}
 \varphi(x) =  \sum_{j=1}^3 (\eta_{j+1} - \eta_{j}) e^{\eta_{j+2} x}  \quad   \mbox{ for } x \in [0, L],  
\end{equation}
and 
\be \label{def-Psi}
\Psi(t, x) =  e^{- i t p}  \varphi(x)  \quad  \mbox{ for } (t, x) \in \mR \times [0, L]. 
\ee
Then $\Psi$ is a solution of the linearized KdV equation thanks to \eqref{eta-p}, i.e., 
\begin{equation} \label{Psi-1}
\Psi_t (t, x) + \Psi_x(t, x) + \Psi_{xxx}(t, x) = 0,   
\end{equation}
and satisfies the boundary conditions 
\begin{equation}\label{Psi-2}
\Psi(t, 0) = \Psi(t, L) = \Psi_x(t, 0) = \Psi_x (t, L) = 0.  
\end{equation}

Multiplying the equation of $y_2$ with $\Psi$, integrating by parts, and using \eqref{Psi-1} and \eqref{Psi-2}, 
one can check that
\be \label{projection-E}
2\int_0^L y_2 (L, x)  \varphi(x) e^{-i p T} \, d x  - 2  \int_0^L y_2 (0, x) \varphi(x) \, dx 
=  \int_0^L \int_{0}^{T} |y_1(t, x)|^2 \varphi_x(x) e^{- i p t} \diff t \diff x.  
\ee
Thus the LHS of \eqref{projection-E} gives the information of the projection of $y_2$ into $\cM_N$. 
The idea in \cite{Cerpa07, CC09} with its root in \cite{CC04}  is then to find $u_1$ and $u_2$ such that, if
$y_1(0, \cdot)=y_2(0, \cdot)=0$, then $y_1(T, \cdot) = 0$ and the $L^2(0, L)$-orthogonal projection of
$y_2(T)$ on $\M_N$ is a given (non-zero) element in $\M_N$ for some time $T>0$. In \cite{CC04}, the authors  needed to make
an expansion up to the order $3$ since $y_2 $ belongs to the orthogonal space of $\M_N$  in
this case. In \cite{CC04, Cerpa07, CC09}, delicate contradiction arguments are involved to capture the structure of the KdV systems.

The analysis in \cite{CKN-24} has the same root. Nevertheless,
instead of using contradiction arguments, our strategy is to characterize all possible
$u_1$ which steers 0 at time 0 to $0$ at time $T$.  This is done by
taking the  Fourier transform with respect to time  of the solution $y_1$  and applying
Paley-Wiener's theorem. With Coron and Koenig \cite{CKN-24}, we showed that 
\be \label{projection}
\int_0^L \int_{0}^{T} |y_1(t, x)|^2 \varphi_x(x) e^{- i p t} \diff t \diff x  = 
\int_{\mR}  \hu_1(z) \overline{\hu_1(z - p )} \int_0^L B(z, x) \diff x \diff z 
\ee
for some function $B(x, z)$ (see, e.g., \Cref{lem-1}) and that, for some complex number $E = E(L)$, it holds  
\be \label{eq-B-1}
\int_0^L B(z, x) \diff x = \frac{E}{|z|^{4/3}} + O(|z|^{-5/3}) \mbox{ for $z \in \mR$ with large $|z|$}. 
\ee
Moreover,  the constant $E$ is not 0 if and only if $2k + l \not \in 3 \N_*$. 
This implies that there are directions in $\M_N$ which cannot be reached via $y_2$ if the time is sufficiently small when $2k + l \not \in 3 \N_*$. Using this observation, we then implement a method to prove
the obstruction to  the local, small-time, null-controllability of the KdV system. The idea is to bring the nonlinear context to the one, based on the power series expansion approach. The proof also relies on  new estimates for the
linear and nonlinear KdV systems using low regularity data. The analysis of this part involves a
connection between the linear KdV equation and the linear KdV-Burgers equation as previously
used by Bona et al. \cite{Bona09} and is inspired by the work of
Bourgain \cite{Bourgain93}, and  Molinet and Ribaud \cite{MR02}. 

\medskip 

The first result of this paper is an improvement of the result of Coron, Koenig, and Nguyen \cite{CKN-24}. 

\begin{theorem}\label{thm-N-1} 
Assume that \eqref{def-L} holds for some 
$k, \,  l \in \N_*$ such that 
\be
2 k + l \not \in  3
\N_*. 
\ee
There exist $T_0>0$, $\eps_0 > 0$, and a real function $\Phi \in C^\infty([0, L]) \cap H_0^1(0, L)$ such
that,  for all $\delta >0$, there is $\eps < \delta$ such that for all $u \in H^{1/2}(0, T_0)$ with $\| u\|_{H^{1/2}(0, T_0)} < \eps_0$, we
have
\[
y(T_0, \cdot) \not \equiv 0,
\]
where $y \in X_T$ is the unique solution of \eqref{sys-KdV-N} with $y_0 = \eps \Phi$ in $[0, L]$. 
\end{theorem}

Here and in what follows, for $T> 0$, denote
\begin{equation}\label{def-XT}
X_T: = C\big([0, T]; L^2(0, L) \big) \cap L^2\big((0, T); H^1(0, L) \big)
\end{equation}
equipped with the corresponding norm: 
$$
\|y\|_{X_T} =\|y \|_{C\big([0, T]; L^2(0, L) \big)} + \| y\|_{L^2\big((0, T); H^1(0, L) \big)}. 
$$

The proof of \Cref{thm-N-1} is given in \Cref{sect-thm-N} using the approach of Coron, Koenig, and Nguyen \cite{CKN-24}. In the analysis, we also incorporate some refinements on the way of bringing the nonlinear context to the one based on the power series expansion approach as implemented in the previous work of the author \cite{Ng-KdV-D}. 

\medskip 

In comparison with the result in \cite{CKN-24}, the controls in \Cref{thm-N-1} are only required to be small with respect to $H^{1/2}$-norm instead of $H^1$-norm, and the compatibility condition at the initial time (see \cite[Theorem 1.2]{CKN-24}) is not imposed in \Cref{thm-N-1}. 

%

\medskip 

Since $E = 0$ in the case $2 k + l  \in 3 \N_*$, to understand the situation when $2 k + l \in 3 \N_*$, one needs to  determine the next-order term in the asymptotic expansion of
$\int_0^L B(z, x) \diff x$ for $z \in \mR$ with large $|z|$. To this end, we show (see \Cref{lem-B}) that, for some non-zero complex number $E_1$, it holds, when $\dim \cM_N \ge 2$,   
\be \label{eq-B-2}
\int_0^L B(z, x) \diff x = \frac{E}{|z|^{4/3}} + \frac{E_1}{|z|^2} + O(|z|^{-7/3}) \mbox{ for $z \in \mR$ with large $|z|$.}
\ee
Using this, we can apply the approach in \cite{CKN-24} to obtain the following result on the case $2k + l \in 3 \N_*$. 

\begin{theorem}\label{thm-N-2} Assume that \eqref{def-L} holds for some $k, l \in \N_*$ with 
\be
2k + l \in 3 \N_*
\ee
There exist $T_0>0$, $\eps_0 > 0$, and a real function $\Phi \in C^\infty[0, L] \cap H_0^1(0, L)$ such
that,  for all $\delta >0$, there is $\eps < \delta$ such that for all $u \in H^{7/6}(0, T_0)$ with $\| u\|_{H^{7/6}(0, T_0)} < \eps_0$, we
have
\[
y(T_0, \cdot) \not \equiv 0,
\]
where $y \in X_T$ is the unique solution of \eqref{sys-KdV-N} with $y_0 = \eps \Phi$ in $[0, L]$. 
\end{theorem}

\begin{remark} \rm  In \Cref{thm-N-2}, we do not impose the compatibility condition at time $0$.  
\end{remark}

The proof of \Cref{thm-N-2} is also included in \Cref{sect-thm-N}, where the proofs of both \Cref{thm-N-1} and \Cref{thm-N-2} are presented.

\medskip 

We next discuss the time used to reach the local exact controllablility of the KdV system \eqref{sys-KdV-N}.  Consider first the case $\dim \cM_N = 2$. The unreachable space $\cM_N$ of the linearized system is then given by 
\be
\mbox{span } \Big\{\Re \varphi (x), \Im \varphi(x) \Big\},   
\ee
where $\varphi$ is defined in \eqref{def-psi}. 

Here and in what follows, for a complex number $z$, $\Re(z)$, $\Im(z)$, and $\bar z$ denote the real part, the imaginary part, and the conjugate of $z$, respectively.

Assume that at a time $\tau > 0$,  there exists a control $u_1$ which brings $y_1$ from 0 to 0, and let $y_2$, starting from $0$ enter into a non-zero direction $\theta$ of $\cM_N$. Recall that $y_1$ and $y_2$ are solutions of the system \eqref{eq:first_order} and \eqref{eq:second_order}, respectively. Then there exists a control $\widetilde u_1$ that again drives $y_1$ from 0 to 0 while steering $y_2$ from 0 into an arbitrary prescribed direction of $\mathcal{M}_N$ at time $\tau + \pi/p$ with $p$ defined by \eqref{def-p}. This is based on the fact that, if $\theta(x) = \Re \varphi \cos (p s_0) + \Im \varphi \sin (p s_0)$ for some $s_0 \in [0, 2 \pi)$ then    
\begin{multline}
\mbox{$ \Re \varphi  \cos (pt + p s_0) + \Im \varphi  \sin (pt + p s_0)$ \mbox{ is a solution of the linearized KdV equation}}\\[6pt]
\mbox{ with zero control and with the initial data $\theta$.}  
\end{multline}
Having this observation, Cerpa \cite{Cerpa07} proved that the control system \eqref{sys-KdV-N} is locally exactly controllable in any time strictly greater than $\pi/ p$ in the case $\dim \cM_N = 2$. Results for $\dim \cM_N \ge 3$ were later derived by Cerpa and Cr\'epeau \cite{CC09} in the same spirit. 

We now recall the time obtained by Cerpa and Cr\'epeau in \cite{CC09} when $\dim \cM_N \ge 2$.  Let $L \in \cM_N$. There exists exactly $n_L \in \N_*$ pairs $(k_m, l_m) \in \N_* \times \N_*$ ($1 \le m \le n_L$) such that $k_m \ge l_m$,  and 
\begin{equation}
L = L (k_m, l_m):= 2 \pi \sqrt{\frac{k_m^2 + k_m l_m + l_m^2}{3}}. 
\end{equation}
For $1 \le m \le n_L$, set 
\begin{equation}\label{def-pm}
p_m = p(k_m, l_m) := \frac{(2k_m + l_m)(k_m - l_m)(2 l_m + k_m)}{3 \sqrt{3}(k_m^2 + k_m l_m + l_m^2)^{3/2}}, 
\end{equation}
 and   denote 
 \begin{equation}
 n_L^{>}= \Big\{p_m \mbox{ given by } \eqref{def-pm} \mbox{ such that } p_m> 0 \Big\}. 
\end{equation}
Then, see, e.g., \cite{CC09} (see also \cite{Cerpa14}),  
$$
n_L - 1 \le n_L^{>} \le n_L, 
$$
and $n_L^{>} = n_{L} - 1$ if and only if the dimension of $\cM_N$ is odd.  In what follows, we consider the convention that 
$$
p_1 > p_2 \dots > p_{n_L^{>}} > 0. 
$$
Set 
\be \label{def-T-CC}
T^{>} = \pi \sum_{m = 1}^{n_L^{>}} \frac{1}{p_m} (n_L^{>} + 1 - m). 
\ee
Cerpa and Cr\'epeau \cite{CC09} proved that the KdV control system \eqref{sys-KdV-N} is locally exactly controllable for initial and final datum in $L^2(0, L)$ using controls in $L^2(0, T)$ for any time larger than $T^{>}$. The intuition behind this result is as follows. 
Define $\varphi_m$ by \eqref{def-psi} with $(k, l) = (k_m, l_m)$ and denote 
\be
M_{m} = \mbox{span } \{\Re \varphi_m, \Im \varphi_m \}.  
\ee
Then 
\be \label{def-MN}
\cM_N = \mbox{span} \Big\{ \varphi; \varphi \in M_m \mbox{ for some $1 \le m \le n_L$} \Big\}. 
\ee
Moreover,
$$
\dim M_m = 2 \mbox{ for } 1 \le m \le n_L^{>}, 
$$
and 
$$
\dim M_m = 1 \mbox{ for } n_L^{>} < m \le n_L. 
$$
The proof of this general case is based on applying to $M_m$ the results in the case $\dim \cM_N = 1$ and in the case $\dim \cM_N = 2$. The time given in \eqref{def-T-CC} is larger than the sum of these times due to a synchronizing process is used, see \cite{CC09}.  

\medskip 
To our knowledge, the time $T^{>}$ given in \eqref{def-T-CC} is the best time obtained so far to reach the local controllability when $\cM_N \ge  2$. It has been asked \cite{CKN-24} if this time is optimal. To our knowledge, the answer is not known even in the case $\dim \cM_N = 2$. 
In this paper, we show that the time $T^{>}$ is not optimal except for two particular lengths. More precisely, we establish the following result. 

\begin{theorem} \label{thm-time}
Let $L \in \cN_N$ with $\dim \cM_N \ge 2$. Assume that 
\be \label{thm-time-as}
\mbox{$L \not \in \big\{L(2, 1), L(3, 1) \big\} = \left\{2 \pi \sqrt{7/3}, 2 \pi \sqrt{13/3} \right\} $. 
}
\ee 
There exists $0< T_0 < T^{>}$ such that the control system \eqref{sys-KdV-N} is locally exactly controllable in time $T_0$ for initial and final datum in $L^2(0, L)$ using controls in $L^2(0, T_0)$.  
\end{theorem}

\Cref{thm-time} is proved in \Cref{sect-thm-time}. The idea is to use \Cref{lem-B} to gain the control time to reach $M_m$ with $\dim M_m = 2$. The key ingredient is given in \Cref{pro-time1}.  The approach of Cerpa and Cr\'epeau \cite{CC09} then allows one to obtain the result.  The requirement \eqref{thm-time-as} seems technical. A refinement of the arguments used here might also yield similar results for these two lengths, but we are not currently able to achieve this.

\subsection{Related results}

A very related control system is the one using the right Dirichlet control: 
\begin{equation}\label{sys-KdV-D}\left\{
\begin{array}{c}
y_t (t, x) + y_x (t, x) + y_{xxx} (t, x) + y (t, x) y_x(t, x) = 0  \mbox{ for } t
\in (0, T), \, x \in (0, L), \\[6pt]
y(t, x=0) =  y(t, x=L) = 0\mbox{ for } t \in (0, T), \\[6pt]
y_x(t, L )  = u \mbox{ for } t \in (0, T), \\[6pt]
y(t = 0, x)  = y_0 (x) \mbox{ for } x \in  (0, L).
\end{array} \right. 
\end{equation}
This control problem was first investigated by Glass and Guerrero \cite{GG10}. To this end, in the spirit of Rosier's work mentioned above, they introduced the corresponding set of critical lengths \footnote{The letter $D$ stands for the Dirichlet boundary control.} 
\begin{equation}\label{def-ND}
\cN_D = \Big\{L \in \mR_+; \, \exists z_1, z_2 \in \mC:  \eqref{cond-z1z2} \mbox{ holds} \Big\}, 
\end{equation}
where 
\begin{equation}\label{cond-z1z2}
z_1 e^{z_1} = z_2 e^{z_2} = - (z_1+ z_2) e^{-(z_1  + z_2)} \quad \mbox{ and } \quad  L^2 = -(z_1^2 + z_1 z_2 + z_2^2). 
\end{equation}
They proved that the set $\cN_D$ is infinite and has no accumulation point. Concerning the control system \eqref{sys-KdV-D},  Glass and Guerrero proved that the corresponding linearized KdV system is small-time, exactly controllable with initial and final data in $H^{-1}(0, L)$ using controls in $L^2(0, T)$  if $L \not \in \cN_D$. Developing this result, they also established that the control system \eqref{sys-KdV-D}  is small-time locally controllable for initial and final data in $L^2(0, L)$ and controls in $H^{1/6_-}$ \footnote{Controls in $H^{1/6_-}(0, T)$ means controls in $H^{1/6 - \eps}(0, T)$ for all $\eps > 0$.}  for {\it non-critical} lengths, i.e., $L \not \in \cN_D$.  

\medskip
The local controllability properties of \eqref{sys-KdV-D} for critical lengths have been obtained very recently by the author \cite{Ng-KdV-D}. We derived a characterization of the critical lengths (see \cite[Proposition 4.1]{Ng-KdV-D}). We also showed  that the unreachable space $\cM_D$ in this case is always of dimension 1 (see \cite[Theorem 1.1]{Ng-KdV-D}), which is quite distinct with the case where the Neumann controls are used. We proved that the control system \eqref{sys-KdV-D} is not locally null controllable in small time for controls in $H^{1/2}(0, T)$ for all critical lengths (see \cite[Theorem 1.2]{Ng-KdV-D}). We also present  a criterion on the local controllability property of \eqref{sys-KdV-D} in finite time (see \cite[Theorem 1.3]{Ng-KdV-D}). As a consequence, we showed that there are critical lengths for which the system is not locally controllable in small time but locally controllable in finite-time (see \cite[Remark 1.3]{Ng-KdV-D}). The analysis has its roots in \cite{CKN-24}. The novelty is on a new look at the corresponding variant of $\int_0^L B(z, x) \diff x$ given in \eqref{eq-B-2}, and various sharp estimates for the solutions of the linearized KdV equation imposed in a bounded interval. The way of bringing the nonlinear context to the one based on the power series expansion approach in \cite{Ng-KdV-D} will be used in the proof of \Cref{thm-N-1,thm-N-2}. 

\medskip 
Recently, with Tran \cite{Ng-Tr-Burgers-25}, we apply the approach of Coron, Koenig, and Nguyen \cite{Coron-Koenig-Nguyen}, and Nguyen \cite{Ng-KdV-D} to study the local null controllability of the Burgers control system $y_t - y_{xx} + y y_x = u(t)$ on a bounded interval imposed by the zero Dirichlet boundary condition. The control $u$ is inside the domain but depends only on $t$. 
It is known from the work of Marbach \cite{Marbach18} that this control system is not locally null controllable in small time using controls in $L^2(0, T)$. In \cite{Ng-Tr-Burgers-25}, we make one step further by proving that the system is not locally null controllable in {\it finite-time} as well using controls in $[H^{3/4}(0, T)]^*$. 

\medskip 

Concerning \Cref{thm-N-2}, in a recent preprint, Niu and Xiang \cite{NX25} claim in the case $2k + l \in 3 \N_*$ that 
the control system \eqref{sys-KdV-N} is not small-time locally null-controllable with controls in $H^{4/3}(0, T)$ under a compatibility condition  (compared with \Cref{thm-N-2} where the controls in $H^{7/6}(0, T)$ are used and no compatibility condition is required).  
Their work closely follows that of Coron, Koenig, and Nguyen \cite{CKN-24} but has some differences, and some analysis in \cite{NX25} must be done more carefully.  We briefly comment on these points.  To 
compute the next-order term in the asymptotic expansion of
$\int_0^L B(z, x) \diff x$ (see \eqref{projection} and \eqref{eq-B-1}) for $z \in \mR$ with large $|z|$ (since $E = 0$ in this case), the authors of \cite{NX25} repeatedly used  \cite[Lemma 3.10]{NX25} (a variant of [CKN24b, Lemma 3.1]). It should be noted, however, that \cite[(3.9)–(3.10) in Lemma 3.10]{NX25} appear to be incompatible; for example, compare them with the computation of $\lambda_2 + \tlambda_2$ in \cite[page 20]{CKN-24} and in \cite[page 30]{NX25} (see also Lemma 3.2 in this paper). To obtain an analogue of \cite[Proposition 3.6]{CKN-24} (see also \Cref{pro-monotone-1} of this paper), the authors of \cite{NX25} employ the theory of pseudo-differential operators (see the proof of \cite[Proposition 4.4]{NX25}). Certain difficulties arise for this approach because the lower-order term in the pseudo-differential expansion must be fully accounted for; therefore, precise estimates that include the dependence on the support of the control are required. 
Consequently, the analysis in \cite{NX25} must be carried out with greater care in this regard.

\medskip 
It is worth mentioning that the power expansion method was also used to study the controllability of the water tank problem \cite{Coron-Koenig-Nguyen-WT}, of the Schr\"odinger equation using bilinear controls \cite{Beauchard05,BL10,BM14,Bournissou23,Gherdaoui25}. It is also used  to study the decay of solutions of KdV systems for critical lengths for which the energy of the solutions of the linearized system is conserved \cite{Ng-decay}.

\subsection{Organization of the paper} 
The paper is organized as follows. In \Cref{sect-0-0}, we recall the characterization of  the controls that steer $0$ to $0$ for the linearized KdV system at a given time from \cite{CKN-24}. In \Cref{sect-dir}, we study
attainable directions for small time via the power series approach. The proofs of \Cref{thm-N-1} and \Cref{thm-N-2} are given in \Cref{sect-thm-N}. The proof of \Cref{thm-time} is finally given in \Cref{sect-thm-time}. A technical lemmas on interpolation inequalities used in the proofs of \Cref{thm-N-1} and \Cref{thm-N-2} is given in the appendix.

\section{Properties of  controls steering \texorpdfstring{$0$ at time $0$ to $0$ at
time $T$}{0 at time 0 to 0 at time T}} \label{sect-0-0}

In this section,  we recall the characterization of  the controls that steer $0$ to $0$ for the linearized KdV system at a given time from \cite{CKN-24}.  Before characterizing controls steering $0$ at time $0$ to $0$ at time $T$, we recall the following definition.

\begin{definition}[\mbox{\cite[Definition 2.2]{CKN-24}}]\label{def:Q-P}
For $z \in \mC$,  let $(\lambda_j)_{1\leq j \leq 3} =  \big(\lambda_j(z) \big)_{1\leq j \leq 3}$ be the three solutions repeated with the multiplicity of
 \begin{equation}\label{eq-lambda}
  \lambda^3 + \lambda + i z = 0.
 \end{equation}
Set
 \begin{equation}\label{eq-defQ}
  Q = Q (z) : =
 \begin{pmatrix}
 1 & 1 & 1\\
 e^{\lambda_1 L } & e^{\lambda_2 L } & e^{\lambda_3 L }\\
 \lambda_1 e^{\lambda_1 L } & \lambda_2 e^{\lambda_2 L } & \lambda_3 e^{\lambda_3 L }
 \end{pmatrix},
\end{equation}
\begin{equation}\label{def-P}
P =  P(z) : =  \sum_{j=1}^3  \lambda_j(e^{\lambda_{j+2} L } - e^{\lambda_{j+1} L})
 = \det
 \begin{pmatrix}
  1&1&1 \\
  e^{\lambda_1L}&e^{\lambda_2L}&e^{\lambda_3L}\\
  \lambda_1& \lambda_2& \lambda_3
 \end{pmatrix},
 \end{equation}
 and
 \begin{equation}\label{def-Xi}
\Xi = \Xi(z) := -  (\lambda_2 - \lambda_1) (\lambda_3 - \lambda_2) (\lambda_1 - \lambda_3) =
\det \begin{pmatrix}1&1&1\\ \lambda_1&\lambda_2&\lambda_3\\
\lambda_1^2&\lambda_2^2&\lambda_3^2\end{pmatrix}, 
\end{equation}
with the convention $\lambda_{j+3} = \lambda_{j}$ for $j \ge 1$. 
\end{definition}

\begin{remark}\label{rk:PG_holomorphic} \rm
The matrix $Q$, and the quantities $P$ and $\Xi$  are antisymmetric  with respect
to $\lambda_j$ ($j=1, 2, 3$), and their definitions depend on a choice of  the order of  $(\lambda_1, \lambda_2, \lambda_3)$. Nevertheless, we later consider a product  of either $P$, $\Xi$,  or $\det Q$ with another antisymmetric function of
$(\lambda_j)$, and these quantities therefore make sense (see, e.g., \eqref{def-hy}, \eqref{def-dxhy}). The definitions of $P$, $\Xi$,  and $Q$ are only  understood in these contexts.
\end{remark}

In what follows, for an appropriate function $v$ defined
on $\mR_+ \times (0, L)$, we extend $v$ by $0$ on $\mR_-\times (0,L)$ and we denote by
$\hat v$ its Fourier transform with respect to $t$, i.e., for $z \in \mC$,
\begin{equation*}
\hat v(z, x) = \frac{1}{\sqrt{2 \pi} }\int_0^{+\infty} v(t, x) e^{- i z t} \diff t.
\end{equation*}

We  have the following result whose proof is based on considering the Fourier transform of the solution in time. 

\begin{lemma}[\mbox{\cite[Lemma 2.4]{CKN-24}}] \label{lem-form-sol} Let $u \in L^2(0, + \infty)$  and let $y \in C
\big([0,
+ \infty); L^2(0, L) \big)$ \\ $ \cap L^2_{\loc}\big( [0, + \infty); H^1(0, L) \big)$ be the
unique solution of
\begin{equation}\label{sys-y}\left\{
\begin{array}{cl}
y_t (t, x) + y_x (t, x) + y_{xxx} (t, x) = 0 &  \mbox{ in } (0, +\infty) \times (0, L),
\\[6pt]
y(t, x=0) = y(t, x=L) = 0 & \mbox{ in }  (0, +\infty), \\[6pt]
y_x(t , x= L) = u(t) & \mbox{ in } (0, +\infty),
\end{array}\right.
\end{equation}
with
\begin{equation}\label{IC-y}
y(t = 0, \cdot) =  0 \mbox{ in } (0, L).
\end{equation}
Then, outside of a discrete set $z\in \mR$, we have
\begin{equation}\label{def-hy}
\hy (z, x) =  \frac{\hu }{\det Q}  \sum_{j=1}^3 (e^{\lambda_{j+2} L } - e^{\lambda_{j+1}
L
}) e^{\lambda_j x}\mbox{ for a.e. } x \in (0, L),
\end{equation}
and in particular,
\begin{equation}\label{def-dxhy}
 \partial_x\hy(z,0)  = \frac{\hu (z) P(z)}{\det Q(z)}.
\end{equation}
\end{lemma}

\begin{remark}\label{rem-form-sol} \rm Assume that $\hu(z, \cdot)$ is well-defined for $z
\in \mC$ (e.g. when $u$ has a compact support). Then the conclusions of
\Cref{lem-form-sol} hold  outside of a discrete set $z\in \mC$.
\end{remark}

As mentioned in \Cref{rk:PG_holomorphic}, the maps $P$ and $\det Q$ are antisymmetric functions with respect to $\lambda_j$. It is hence convenient to consider $\partial_x
\hy(z, 0)$ under the form
\begin{equation}\label{eq:hyx_meromorphic}
\partial_x \hy(z, 0) = \frac{\hu (z) G(z) }{H(z)},
\end{equation}
where, with $\Xi$ defined in \eqref{def-Xi},
\begin{equation}\label{def-GH}
G(z) =  P(z)/\Xi(z) \quad \mbox{ and } \quad H(z) = \det Q(z)/\Xi(z), 
\end{equation}
which are entire functions, see \cite[Lemma 2.6]{CKN-24}. 

\medskip

We are ready to recall the characterization of the controls
steering $0$ to $0$ from \cite{CKN-24}. 

\begin{proposition}[\mbox{\cite[Proposition 2.8]{CKN-24}}] \label{pro-Gen} Let $L>0$, $T>0$, and $u \in L^2(0, + \infty)$. Assume that   $u$ has a compact support in $[0, T]$,  and $u$  steers $0$ at the time $0$ to $0$ at the time $T$, i.e., the unique solution $y$ of \eqref{sys-y} and \eqref{IC-y} satisfies $y(T, \cdot) = 0$ in $(0, L)$.  Then $\hu$ and  $\hu G/ H$  satisfy the assumptions of Paley-Wiener's theorem concerning the support in $[-T, T]$, i.e.,
\[
\hu \mbox{ and } \hu G/H \mbox{ are entire functions},
\]
and
\[
|\hu(z)|  + \left|\frac{\hu G(z)}{H(z)} \right| \le C e^{T| \Im(z)|},
\]
for some positive constant $C$.
\end{proposition}

\section{Attainable directions in small time}\label{sect-dir}

Let  $u \in L^2(0,  + \infty)$ and denote $y$ the corresponding solution
of the linear KdV equation~\eqref{sys-y} with zero initial condition. We assume
that $y$ satisfies $y(t, \cdot) = 0$ in $(0, L)$ for $t \ge T$. We have, by
\Cref{lem-form-sol} (and also \Cref{rem-form-sol}), for $z \in \mC$ outside  a discrete
set,
\begin{equation}\label{def-y}
\hat y(z, x) = \hat u(z) \frac{\sum_{j=1}^3(e^{\lambda_{j+1} L } - e^{\lambda_{j} L })e^{\lambda_{j+2} x} }{\sum_{j=1}^3 (\lambda_{j+1} - \lambda_{j}) e^{-\lambda_{j+2} L }}.
\end{equation}
Recall that  $\lambda_j = \lambda_j(z)$ for $j=1, \, 2, \, 3$  are the three solutions of
the equation
\begin{equation}\label{eq-lambda-j}
\lambda^3 + \lambda + i z = 0.
\end{equation}
Let $\eta_1, \, \eta_2, \,  \eta_3 \in i \mR$, i.e.,  $\eta_j \in \mC$ with $\Re(\eta_j) = 0$ for $j=1, \, 2, \, 3$. Define
\begin{equation}\label{def-varphi}
\varphi(x) = \sum_{j=1}^3 (\eta_{j+1} - \eta_{j}) e^{\eta_{j+2} x} \mbox{ for } x \in [0, L],
\end{equation}
with the convention $\eta_{j+3} = \eta_j$ for $j \ge 1$. The following assumption on $\eta_j$ is used repeatedly  throughout the paper:
\begin{equation}\label{pro-eta}
e^{\eta_1  L} = e^{\eta_2 L} = e^{\eta_3 L},
\end{equation}
which is equivalent to $\eta_3 - \eta_2, \eta_2 - \eta_1 \in \dsp  \frac{2 \pi i}{L} \mZ$. 

\medskip 

We have the following result.

\begin{lemma}[\mbox{\cite[Lemma 3.1]{CKN-24}} 
]\label{lem-1} Let   $p \in \mR$ and let $\varphi$ be defined by \eqref{def-varphi}.  Set, for $(z, x) \in  \mC \times [0, L]$,
\begin{equation}\label{def-B}
B(z, x) = \frac{\sum_{j=1}^3(e^{\lambda_{j+1} L } - e^{\lambda_{j} L })e^{\lambda_{j+2} x} }{\sum_{j=1}^3 (\lambda_{j+1} - \lambda_{j}) e^{-\lambda_{j+2} L }} \cdot  \frac{\sum_{j=1}^3(e^{\tlambda_{j+1} L } - e^{\tlambda_{j} L })e^{\tlambda_{j+2} x} }{\sum_{j=1}^3 (\tlambda_{j+1} - \tlambda_{j}) e^{-\tlambda_{j+2} L }}  \cdot \varphi_x(x),
\end{equation}
where $\tlambda_j = \tlambda_j (z)$ ($j=1, \, 2, \, 3$) denotes the conjugate of the roots of \eqref{eq-lambda-j}  with $z$ replaced by $z - p$ and with the use of  the convention $\tlambda_{j+3} = \tlambda_j$ for $j \ge 1$.
Let $u \in L^2(0, + \infty)$ and let $y \in C \big([0, + \infty); L^2(0, L) \big) \cap L^2_{\loc}\big( [0, + \infty); H^1(0, L) \big)$ be the unique solution of \eqref{sys-y} and \eqref{IC-y}.  Assume that $u(t) = 0$  and  $y(t, \cdot) = 0$ for large $t$. 
Then
\begin{equation}\label{identity-1}
\int_0^L \int_{0}^{+ \infty} |y(t, x)|^2 \varphi_x(x) e^{- i p t} \diff t \diff x  =
\int_{\mR}  \hu(z) \overline{\hu(z - p )} \int_0^L B(z, x) \diff x \diff z .
\end{equation}
\end{lemma}

We next state  the behaviors of $\lambda_j$ and $ \tlambda_j$ given in \Cref{lem-1} for $z \in \mR$ with large $|z|$. 
These asymptotics are a direct consequence of the equation~\eqref{eq-lambda} satisfied by
the $\lambda_j$. A variant of this was given in \cite[Lemma 3.3]{CKN-24}.

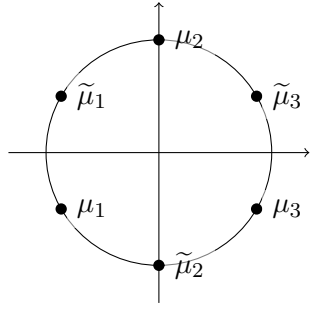
\begin{figure}[htp]
 
\begin{tikzpicture}[scale=0.5]
  \draw[->] (-4,0) -- (4,0);
  \draw[->] (0,-4) -- (0,4);
  
  \draw (0,0) circle[radius=3];
  
  \foreach \j in {1,2,3} {
   \fill (-30-120*\j:3) circle[fill, radius=0.15] 
     ++(0.15,0) node[right, fill=white, fill opacity=0.5, text opacity = 1]{$\mu_\j$};
   \fill (30+120*\j:3) circle[fill, radius=0.15] 
     ++(0.15,0) node[right, fill=white, fill opacity=0.5, text opacity = 1]{$\tmu_\j$}; 
     }
\end{tikzpicture}

 \caption{The roots $\lambda_j$ of $\lambda^3 + \lambda + i z = 0$ satisfy, when
$z>0$ is large, $\lambda_j \sim \mu_j z^{1/3}$ where $\mu_j^3 = -i$.}
\label{fig:mu}
\end{figure}

\begin{lemma}\label{th:lambda_asym}
For  $p \in \mC$ and $z > 0$ large, let $\lambda_j$ and $\tlambda_j$ with $j=1, \, 2, \, 3$
be given in  \Cref{lem-1}.
Consider the convention $\Re(\lambda_1) < \Re(\lambda_2) < \Re(\lambda_3)$ and similarly
for $\tlambda_j$. We have,  in the
limit $z \to + \infty$,
\begin{equation}
\lambda_j = \mu_jz^{1/3} - \frac1{3\mu_j}z^{-1/3} + O(z^{-5/3}) \quad \text{ with } \mu_j
= e^{-i\pi/6-2ji\pi/3},
\end{equation}
\begin{equation}
\tlambda_j = \tmu_j z^{1/3} - \frac1{3\tmu_j} z^{-1/3} -  \frac{1}{3} \tmu_j p z^{-2/3}  - \frac{1}{9 \tmu_j} p z^{-4/3} +O(z^{- 5/3})  \quad  \text{ with }
\tmu_j = e^{i\pi/6+2ij\pi/3}
\end{equation}
(see Figure~\ref{fig:mu} for the geometry of $\mu_j$ and $\tmu_j$). Here $z^{1/3}$ denotes the real cube root of $z$.  
\end{lemma}

\begin{remark} \rm In \Cref{lem-1}, we only defined $\tlambda_j$ for $p \in \mR$. The same definition is also used for $p \in \mC$. 
\end{remark}

\begin{remark} \rm The case $z \in \mR$ with large $|z|$ can be obtained from the $z > 0$ with large $|z|$ by 
considering the conjugate of suitable quantities.  
\end{remark}

%
%

We are ready to establish the behavior of
\[
\int_0^L B(z, x) \diff x
\]
for $z \in \mR$ with  large $|z|$, which is one of the main ingredients for the analysis in this section.

\begin{lemma} \label{lem-B}  Let   $p \in \mR \setminus \{0 \}$,  and let $\varphi$ be defined by \eqref{def-varphi} where 
\be \label{lem-B-p1}
\eta_j^3 + \eta_j - ip = 0 \mbox{ for } 1 \le j \le 3.  
\ee
Assume that $\Re{\eta_j} = 0$ and  \eqref{pro-eta} holds. Let $B$ be defined by \eqref{def-B}. We have
\be \label{lem-B-cl1}
 \int_0^L B(z, x) \diff x  = E z^{-4/3}  + E_1 |z|^{-2} 
+ O(|z|^{-7/3})
\mbox{ for $z\in \mR$ with large $|z|$},
\ee
where 
\be
E = \frac{1}{3} (e^{\eta_1 L} - 1)\left( - \frac{2}{3} \Gamma  -  \frac{1}{3} \Lambda   \right) 
\ee
and
\be
E_1 = - \frac{1}{3} (1 + i p L) E + F. 
\ee
Here 
\be
\Gamma = \sum_{j=1}^3 (\eta_{j+1} - \eta_j) \eta_{j+2}^2,  \quad \Lambda = i p \sum_{j=1}^3 \frac{\eta_{j+1} - \eta_j}{\eta_{j+2}}, 
\ee
and 
\be
F =  \frac{2}{27} (\Lambda - \Gamma) (e^{\eta_1 L} - 1)  +  \frac{1}{9} i p Le^{\eta_1 L}  \left( - \frac{2}{3} \Gamma  -  \frac{1}{3} \Lambda   \right).
\ee
Moreover, if $e^{\eta_1 L} = 1$, it holds 
\be \label{lem-B-cl2}
 \int_0^L B(z, x) \diff x  = F |z|^{-2}  + F_1 |z|^{-8/3} 
+ O(|z|^{-9/3})
\mbox{ for $z\in \mR$ with large $|z|$},
\ee
where 
\be
F_1 = F  \left(- \frac{2}{3}  - \frac{ip L}{3} + \frac{2 (\Gamma - \Lambda) }{3 (2 \Gamma + \Lambda)} \right). 
\ee
\end{lemma}

\begin{remark} \rm Note that if $e^{\eta_1 L } = 1$, then 
$$
E = 0 \quad \mbox{ and } \quad E_1 = F  = - \frac{1}{27} ip L(2 \Gamma + \Lambda).
$$
\end{remark}

\begin{remark} \rm
Since 
\be  \label{lem-B-p2}
\sum_{j=1}^3 (\eta_{j+1} - \eta_j) \eta_{j+2} = 0,
\ee
it follows from \eqref{lem-B-p1} that 
\be  \label{lem-B-p12}
\sum_{j=1}^3 (\eta_{j+1} - \eta_j) \eta_{j+2}^3 = 0,  \quad \sum_{j=1}^3 (\eta_{j+1} - \eta_j) \eta_{j+2}^4 = - \sum_{j=1}^3 (\eta_{j+1} - \eta_j) \eta_{j+2}^2 = - \Gamma. 
\ee
\end{remark}

\begin{proof}
We first deal with the case where $z$ is positive and large. We use the convention in \Cref{th:lambda_asym} for $\lambda_j$ and $\tlambda_j$.
Consider the denumerator of $B(z, x)$. We have, by \Cref{th:lambda_asym},
\begin{multline}\label{den-B}
\frac{1}{\sum_{j=1}^3 (\lambda_{j+1} - \lambda_{j}) e^{-\lambda_{j+2} L }} \cdot
\frac{1}{ \sum_{j=1}^3 (\tlambda_{j+1} - \tlambda_{j}) e^{-\tlambda_{j+2} L }} \\[6pt]
= \frac{e^{\lambda_1 L } e^{\tlambda_1 L}}{(\lambda_3 - \lambda_2)  (\tlambda_3 -
\tlambda_2)  } \Big( 1 + O \big(e^{-C |z|^{1/3}} \big) \Big).
\end{multline}

We next deal with the numerator of $B(z, x)$. Set, for $(z, x) \in \mR \times (0, L)$,
\begin{equation}\label{def-fg}
f(z, x) = \sum_{j=1}^3(e^{\lambda_{j+1} L } - e^{\lambda_{j} L })e^{\lambda_{j+2} x},
\quad g (z, x) = \sum_{j=1}^3(e^{\tlambda_{j+1} L } - e^{\tlambda_{j} L
})e^{\tlambda_{j+2} x},
\end{equation}
\footnote{The index $m$ stands the
main part.}
\[
f_m(z, x) =  - e^{\lambda_3 L} e^{\lambda_2 x} + e^{\lambda_2 L}  e^{\lambda_3 x} +
e^{\lambda_3 L } e^{\lambda_1 x}, \quad
g_m (z, x) =  - e^{\tlambda_3 L} e^{\tlambda_2 x} + e^{\tlambda_2 L}  e^{\tlambda_3 x}
+ e^{\tlambda_3 L } e^{\tlambda_1 x}.
\]
We have
\begin{multline*}
 \int_0^L f(z, x) g(z, x) \varphi_x(x) \diff x
 = \int_0^L f_m(z, x) g_m(z, x) \varphi_x(x) \diff x + \int_0^L (f- f_m)(z, x) g_m(z, x)
\varphi_x(x) \diff x \\[6pt]  + \int_0^L f_m(z, x) (g -g_m) (z, x) \varphi_x(x) \diff x +
\int_0^L (f- f_m)(z, x) (g -g_m) (z, x) \varphi_x(x) \diff x.
\end{multline*}
It is clear from \Cref{th:lambda_asym} that
\begin{multline}\label{B-p0}
 \int_0^L |(f- f_m)(z, x) g_m(z, x) \varphi_x(x)| \diff x + \int_0^L |(f- f_m)(z, x) (g
-g_m) (z, x) \varphi_x(x)| \diff x  \\[6pt]  + \int_0^L |f_m(z, x) (g -g_m) (z, x)
\varphi_x(x)| \diff x \le C |e^{(\lambda_3 + \tlambda_3) L }| e^{- C |z|^{1/3}}.
\end{multline}

We next estimate
\begin{equation}\label{eq:int_fm_gm}
\int_0^L f_m(x, z) g_m(x, z) \varphi_x(x) = \int_0^L f_m(x, z) g_m(x, z) \left(
\sum_{j=1}^3 \eta_{j+2} (\eta_{j+1} - \eta_{j}) e^{\eta_{j+2} x} \right) \diff x.
\end{equation}
We first have, by  \eqref{pro-eta} and \Cref{th:lambda_asym},
\begin{multline}\label{B-p1}
\int_0^L \Big(  - e^{\lambda_3 L} e^{\lambda_2 x}  e^{\tlambda_2 L}  e^{\tlambda_3 x}
- e^{\lambda_2 L}  e^{\lambda_3 x}e^{\tlambda_3 L} e^{\tlambda_2 x}
+ e^{\lambda_2 L}  e^{\lambda_3 x}e^{\tlambda_2 L}  e^{\tlambda_3 x} \Big) \\[6pt]
\times \left( \sum_{j=1}^3 \eta_{j+2} (\eta_{j+1} - \eta_{j}) e^{\eta_{j+2} x} \right)
\diff x
=  e^{ (\lambda_3 +   \tlambda_3 + \lambda_2 +  \tlambda_2) L } \Big(e^{\eta_1 L
}T_1(z) + O\big(e^{-C|z|^{1/3}} \big) \Big),
\end{multline}
where
\begin{equation}\label{def-T1}
T_1(z) : =  \sum_{j=1}^3 \eta_{j+2}(\eta_{j+1}-\eta_{j}) \left( \frac{1}{  \lambda_3
+   \tlambda_3   + \eta_{j+2}} - \frac{1}{ \lambda_3  +   \tlambda_2   + \eta_{j+2}} -
\frac{1}{ \lambda_2  + \tlambda_3  + \eta_{j+2}}  \right).
\end{equation}

Let us now deal with the terms of~\eqref{eq:int_fm_gm} that contain both $e^{\lambda_3 L + \tlambda_3 L}$ and (either $e^{\lambda_1 x}$ or $e^{\tlambda_1 x}$).  We obtain, by  \eqref{pro-eta} and \Cref{th:lambda_asym},
\begin{multline}\label{B-p2}
\int_0^L \Big( e^{\lambda_3 L } e^{\lambda_1 x}e^{\tlambda_3 L } e^{\tlambda_1 x}
- e^{\lambda_3 L } e^{\lambda_1 x} e^{\tlambda_3 L} e^{\tlambda_2 x} - e^{\lambda_3 L}
e^{\lambda_2 x} e^{\tlambda_3 L } e^{\tlambda_1 x} \Big) \\[6pt]
\times \left( \sum_{j=1}^3 \eta_{j+2} (\eta_{j+1} - \eta_{j}) e^{\eta_{j+2} x} \right)
\diff x
=  e^{(\lambda_3 + \tlambda_3) L }\Big( T_2(z) + O(e^{-C|z|^{1/3}}) \Big),
\end{multline}
where
\begin{equation}\label{def-T2}
T_2 (z): = \sum_{j=1}^3 \eta_{j+2}(\eta_{j+1}-\eta_{j}) \left( - \frac{1}{\lambda_1
+  \tlambda_1  + \eta_{j+2}} + \frac{1}{\lambda_1   +  \tlambda_2   + \eta_{j+2}} +
\frac{1}{ \lambda_2  +   \tlambda_1   + \eta_{j+2}}\right).
\end{equation}
We  have, by \eqref{pro-eta},
\begin{equation}\label{B-p3}
\int_0^L e^{\lambda_3 L} e^{\lambda_2 x} e^{\tlambda_3 L} e^{\tlambda_2 x}  \left(
\sum_{j=1}^3 \eta_{j+2} (\eta_{j+1} - \eta_{j}) e^{\eta_{j+2} x} \right) \diff x
=  e^{(\lambda_3 + \tlambda_3) L } T_3(z),
\end{equation}
where
\begin{equation}\label{def-T3}
T_3 (z): = \Big( e^{\lambda_2 L  + \tlambda_2 L  + \eta_{1} L } - 1\Big)  \hat T_3 (z) \quad \mbox{ with } \quad 
\hat T_3 (z) =  \sum_{j=1}^3
\frac{\eta_{j+2}(\eta_{j+1}-\eta_{j})  }{\lambda_2   + \tlambda_2   + \eta_{j+2}}.
\end{equation}
The other terms of~\eqref{eq:int_fm_gm} are negligible, because we have
\begin{multline}\label{B-p3-1}
\left| \int_0^L \Big(  e^{\lambda_3 L } e^{\lambda_1 x} e^{\tlambda_2 L} e^{\tlambda_3
x} + e^{\lambda_2 L} e^{\lambda_3 x} e^{\tlambda_3 L} e^{\tlambda_1 x}  \Big) \Big(
\sum_{j=1}^3 \eta_{j+2} (\eta_{j+1} - \eta_{j}) e^{\eta_{j+2} x} \Big) \diff x \right|
\\[6pt]
= |e^{(\lambda_3 + \tlambda_3)L} | O(e^{-C z^{1/3}}).
\end{multline}

We thus need to estimate, after taking $e^{(\lambda_3 + \tlambda_3)L}$ out,  
\begin{multline}
\Big(e^{\lambda_2 L + \tlambda_2 L + \eta_1 L} - 1\Big) T_1(z) + \Big(T_1(z) + T_2(z) \Big) +  \Big( e^{\lambda_2 L + \tlambda_2 L + \eta_1 L} -1 \Big) \hat T_3 (z) \\[6pt]
= (T_1 (z) + T_2(z)) +  \Big(e^{\lambda_2 L + \tlambda_2 L + \eta_1 L} - 1\Big) \Big(T_1 (z) + \hat T_3 (z) \Big).
\end{multline}

%
%

We first deal with $T_1(z)$ given in
\eqref{def-T1}. Using \eqref{lem-B-p2}, we obtain
\begin{align*}
T_1 (z)=  & \sum_{j=1}^3 \eta_{j+2}(\eta_{j+1}-\eta_{j}) \left( \frac{1}{\lambda_3 +
\tlambda_3 + \eta_{j+2}} -  \frac{1}{\lambda_3 + \tlambda_3 } \right) \\[6pt]
& +
 \sum_{j=1}^3 \eta_{j+2}(\eta_{j+1}-\eta_{j})  \left(  - \frac{1}{\lambda_3 +
\tlambda_2 + \eta_{j+2}}
 +  \frac{1}{\lambda_3 + \tlambda_2}  \right) \\[6pt]
&  + \sum_{j=1}^3 \eta_{j+2}(\eta_{j+1}-\eta_{j})  \left(- \frac{1}{\lambda_2 +
\tlambda_3 + \eta_{j+2}}  + \frac{1}{\lambda_2 + \tlambda_3}  \right).
\end{align*}
Using \Cref{th:lambda_asym}, we get
\be \label{T1-p0-0}
T_1 (z)  = \sum_{l=2}^6   (-1)^l a_l T_{1, l} (z) + O(z^{-7/3}), 
\ee
where 
\be \label{def-al}
a_l = \sum_{j=1}^3 \eta_{j+2}^l (\eta_{j+1} -\eta_{j}) \quad \mbox{ and } \quad 
T_{1, l}  (z) =  - \frac{1}{(\lambda_3
+ \tlambda_3)^l} + \frac{1}{(\lambda_3 + \tlambda_2)^l} + \frac{1}{(\lambda_2 +
\tlambda_3)^l} . 
\ee
Note that 
\be
a_3 \mathop{=}^{\eqref{lem-B-p12}} 0. 
\ee
and 
\be
T_{1, 5} (z) \mathop{=}^{\Cref{th:lambda_asym}}  z^{-5/3} \left(- \frac{1}{(\mu_3
+ \tmu_3)^5} + \frac{1}{(\mu_3 + \tmu_2)^5} + \frac{1}{(\mu_2 +
\tmu_3)^5} \right)  + O(z^{-7/3}) = O(z^{-7/3})
\ee
since 
$$
- \frac{1}{(\mu_3
+ \tmu_3)^5} + \frac{1}{(\mu_3 + \tmu_2)^5} + \frac{1}{(\mu_2 +
\tmu_3)^5} = 0. 
$$
It follows from \eqref{T1-p0-0} that 
\be \label{T1-p0}
T_1(z) = a_2 T_{1, 2} (z) + a_4 T_{1, 4} (z)+ a_6 T_{1, 6} (z) + O(z^{-7/3}). 
\ee

From \Cref{th:lambda_asym}, we have 
\begin{multline}\label{lem-B-T1-p2}
T_{1, 2} (z) =  - \frac{1}{(\lambda_3 + \tlambda_3)^2} + \frac{1}{(\lambda_3 + \tlambda_2)^2} +
\frac{1}{(\lambda_2 + \tlambda_3)^2} \\[6pt]
 = z^{-2/3} \Big(- (\mu_3+\tmu_3)^{-2} + (\mu_3+\tmu_2)^{-2} +(\mu_2+\tmu_3)^{-2}
\Big) \\[6pt]
+  \frac{2}{3} z^{-4/3} \left( - \frac{(\mu_3+\tmu_3)^{-2}}{\mu_3 \tmu_3} + \frac{(\mu_3+\tmu_2)^{-2}}{\mu_3 \tmu_2} + \frac{(\mu_2+\tmu_3)^{-2}}{\mu_2 \tmu_3} \right)
\\[6pt]
+ \frac{2}{3}  p z^{-5/3}   \Big( - \tmu_3(\mu_3+\tmu_3)^{-3} +  \tmu_2 (\mu_3+\tmu_2)^{-3} + \tmu_3(\mu_2+\tmu_3)^{-3} \Big) \\[6pt]
+ \frac{1}{3} z^{-2} \left( - \frac{(\mu_3+\tmu_3)^{-2}}{\mu_3^2 \tmu_3^2} + \frac{(\mu_3+\tmu_2)^{-2}}{\mu_3^2 \tmu_2^2} + \frac{(\mu_2+\tmu_3)^{-2}}{\mu_2^2 \tmu_3^2} \right) \\[6pt]
= -  \frac23z^{-2/3} - \frac{4}{9} z^{-4/3}  -\frac{2}{9}  p z^{-5/3} + \frac{1}{9} z^{-2}+ O(z^{-7/3}),
\end{multline}
and
\begin{multline} \label{lem-B-T1-p3}
 T_{1, 4} (z) = - \frac{1}{(\lambda_3 + \tlambda_3)^4} + \frac{1}{(\lambda_3 + \tlambda_2)^4} +
\frac{1}{(\lambda_2 + \tlambda_3)^4} \\[6pt]
 = z^{-4/3} \Big(- (\mu_3+\tmu_3)^{-4} + (\mu_3+\tmu_2)^{-4} +(\mu_2+\tmu_3)^{-4}
\Big) + O(z^{-5/3})
= -  \frac{2}{9}z^{-4/3} + O(z^{-5/3}). 
\end{multline}
Combining  \eqref{T1-p0},  \eqref{lem-B-T1-p2}, and \eqref{lem-B-T1-p3},  and using the fact that $a_2 = \Gamma$ and $a_4 = - \Gamma$ by \eqref{lem-B-p12}, we obtain  
\be\label{T1-final}
T_1   = - \frac{2}{3} \Gamma z^{-2/3}   
- \frac{2}{9} \Gamma z^{-4/3}+ O(z^{-5/3}).
\ee

We next consider  $T_2(z)$ given in \eqref{def-T2}. We have, by \eqref{lem-B-p2},
\begin{align*}
T_2 (z)=  & \sum_{j=1}^3 \eta_{j+2}(\eta_{j+1}-\eta_{j}) \left( -  \frac{1}{\lambda_1
+ \tlambda_1 + \eta_{j+2}} +  \frac{1}{\lambda_1 + \tlambda_1 } \right) \\[6pt]
& +
 \sum_{j=1}^3 \eta_{j+2}(\eta_{j+1}-\eta_{j})  \left(   \frac{1}{\lambda_1 +
\tlambda_2 + \eta_{j+2}} -
   \frac{1}{\lambda_1 + \tlambda_2}  \right) \\[6pt]
&  + \sum_{j=1}^3 \eta_{j+2}(\eta_{j+1}-\eta_{j})  \left( \frac{1}{\lambda_2 +
\tlambda_1 + \eta_{j+2}}  - \frac{1}{\lambda_2 + \tlambda_1}  \right).
\end{align*}
Using \Cref{th:lambda_asym}, we get
\be \label{T2-p0-0}
T_2  (z) = \sum_{l=2}^6  (-1)^l a_l T_{2, l} (z) + O(z^{-7/2}), 
\ee
where $a_l$ is defined in \eqref{def-al} for $2 \le l \le 6$, and 
\be
T_{2, l} (z)  = \frac{1}{(\lambda_1
+ \tlambda_1)^l} - \frac{1}{(\lambda_1 + \tlambda_2)^l} - \frac{1}{(\lambda_2 +
\tlambda_1)^l}. 
\ee
Note that 
\be
a_3 \mathop{=}^{\eqref{lem-B-p12}} 0. 
\ee
and 
\be
T_{2, 5} (z) =  z^{-5/3} \left(\frac{1}{(\mu_1
+ \tmu_1)^5} - \frac{1}{(\mu_1 + \tmu_2)^5} - \frac{1}{(\mu_2 +
\tmu_1)^5} \right)  + O(z^{-7/3}) = O(z^{-7/3}),  
\ee
since 
$$
\frac{1}{(\mu_1
+ \tmu_1)^5} - \frac{1}{(\mu_1 + \tmu_2)^5} - \frac{1}{(\mu_2 +
\tmu_1)^5}  = 0.   
$$
It follows from \eqref{T2-p0-0} that 
\be \label{T2-p0}
T_2 (z) = a_2 T_{2, 2} (z) + a_4 T_{2, 4} (z) + a_6 T_{2, 6} (z) + O(z^{-7/3}). 
\ee

From \Cref{th:lambda_asym}, similar to \eqref{lem-B-T1-p2}, we have 
\begin{multline}\label{lem-B-T2-p2}
T_{2, 2}  (z) =   \frac{1}{(\lambda_1 + \tlambda_1)^2} - \frac{1}{(\lambda_1 + \tlambda_2)^2} -
\frac{1}{(\lambda_2 + \tlambda_1)^2} \\[6pt]
 = z^{-2/3} \Big( (\mu_1+\tmu_1)^{-2} - (\mu_1+\tmu_2)^{-2}  - (\mu_2+\tmu_1)^{-2}
\Big) \\[6pt]
+  \frac{2}{3} z^{-4/3} \left(  \frac{(\mu_1+\tmu_1)^{-2}}{\mu_1 \tmu_1} - \frac{(\mu_1+\tmu_2)^{-2}}{\mu_1 \tmu_2} - \frac{(\mu_2+\tmu_1)^{-2}}{\mu_2 \tmu_1} \right)
\\[6pt]
+ \frac{2}{3}  p z^{-5/3}   \left( \tmu_1(\mu_1+\tmu_1)^{-3} -  \tmu_2 (\mu_1+\tmu_2)^{-3} - \tmu_1(\mu_2+\tmu_1)^{-3} \right)\\[6pt]
+ \frac{1}{3} z^{-2} \left( \frac{(\mu_1+\tmu_1)^{-2}}{\mu_1^2 \tmu_1^2} - \frac{(\mu_1+\tmu_2)^{-2}}{\mu_1^2 \tmu_2^2} - \frac{(\mu_2+\tmu_1)^{-2}}{\mu_2^2 \tmu_1^2} \right) \\[6pt]
=   \frac23z^{-2/3} + \frac{4}{9} z^{-4/3}  + \frac{2}{9}  p z^{-5/3} - \frac{1}{9} z^{-2}+ O(z^{-7/3})
\end{multline}

Since, by \Cref{th:lambda_asym},  
\begin{align*}
\lambda_3 + \tlambda_3  = - (\lambda_1 + \tlambda_1) 
+ O(z^{-2/3}), \\[6pt]
 \lambda_3 + \tlambda_2 = - (\tlambda_1 + \lambda_2)    + O(z^{-2/3}), \\[6pt]
\lambda_2 + \tlambda_3 = - (\tlambda_2 + \lambda_1) + O(z^{-2/3}),  
\end{align*}
it follows that 
\be \label{lem-B-T12-46}
T_{1, 4}(z) + T_{2, 4} (z) = O(z^{-7/3}) \quad \mbox{ and } \quad T_{1, 6} (z) + T_{2, 6}(z) = O(z^{-7/3}). 
\ee
Combing \eqref{T1-p0} and \eqref{T2-p0}, and using \eqref{lem-B-T1-p2}, \eqref{lem-B-T2-p2}, and \eqref{lem-B-T12-46}, we obtain  
\be \label{T1-T2}
T_{1} (z) + T_{2} (z) =    O(z^{-7/3}).
\ee

We finally  consider  $\hat T_3(z)$ given in \eqref{def-T3}. We have, by \Cref{th:lambda_asym}, 
\be \label{lambda2-tlambda2}
\lambda_2 + \tlambda_2 = \frac{1}{3} i p  z^{-2/3}   - \frac{1}{9} i p  z^{-4/3}+ O(z^{-5/3}).
\ee
From \eqref{lambda2-tlambda2}, we derive that 
\begin{multline}\label{T3-p2}
\hat T_3 (z) = \sum_{j=1}^3 \frac{\eta_{j+2}(\eta_{j+1}-\eta_{j})  }{\lambda_2 + \tlambda_2  +
\eta_{j+2}}  =  \sum_{j=1}^3 \frac{\eta_{j+2}(\eta_{j+1}-\eta_{j})  }{\frac{1}{3}i p
z^{-2/3}  - \frac{ip}{9} z^{-4/3} + \eta_{j+2}}  + O(|z|^{-5/3}) \\[6pt]
=   \sum_{j=1}^3 (\eta_{j+1}-\eta_{j})  \left( 1 -  \frac{ip z^{-2/3}}{ 3 \eta_{j+2} } +   \frac{ip}{9 \eta_{j+2}} z^{-4/3} - \frac{p^2}{9 \eta_{j+2}^2} z^{-4/3}
\right) + O(|z|^{-5/3}). 
\end{multline}
Since 
$$
\frac{1}{\eta_{j+2}} + \frac{ip}{\eta_{j+2}^2}  \mathop{=}^{\eqref{lem-B-p1}} \frac{1}{\eta_{j+2}} + \frac{\eta_{j+2}^3 + \eta_{j+2}}{\eta_{j+2}^2}  = \eta_{j+2} + \frac{2}{\eta_{j+2}}, 
$$
it follows from \eqref{lem-B-p2} that 
\begin{multline}
\label{T3-final}
\hat T_3 (z) = 
\sum_{j=1}^3 \frac{\eta_{j+2}(\eta_{j+1}-\eta_{j})  }{\lambda_2 + \tlambda_2  +
\eta_{j+2}}  \\[6pt]
=   - \frac{ i p }{3}   \sum_{j=1}^3 \frac{  \eta_{j+1}-\eta_{j}} {
\eta_{j+2} } z^{-2/3}   + \frac{2 i p}{9}\sum_{j=1}^3 \frac{  \eta_{j+1}-\eta_{j}} {
\eta_{j+2} } z^{-4/3} +  O(z^{-5/3}) = - \frac{1}{3} \Lambda z^{-2/3} + O(z^{-5/3}). 
\end{multline}

Combining \eqref{T1-final} and \eqref{T3-final},  
we thus obtain 
\be \label{T1-T3}
T_1 (z) + \hat T_3 (z) = \left( - \frac{2}{3} \Gamma  -  \frac{1}{3} \Lambda   \right) z^{-2/3}   +  \frac{2}{9} (\Lambda - \Gamma) z^{-4/3} + O(z^{-5/3}).  
\ee
Since
\be \label{lem-T-sum}
\Big(e^{\lambda_2 L + \tlambda_2 L + \eta_1 L} - 1\Big) (T_1 + \hat T_3)\\[6pt]
=   \Big(e^{\eta_1 L} - 1\Big) (T_1 + \hat T_3) + e^{\eta_1 L}\Big(e^{\lambda_2 L + \tlambda_2 L } - 1 \Big) (T_1 + \hat T_3),  
\ee
using \eqref{lambda2-tlambda2}, we derive from \eqref{T1-T2} and \eqref{T1-T3} that
\begin{multline}\label{sum1}
(T_1 + T_2)  + \Big(e^{\lambda_2 L + \tlambda_2 L + \eta_1 L} - 1\Big) (T_1 + \hat T_3)
=  (e^{\eta_1 L} - 1)\left( - \frac{2}{3} \Gamma  -  \frac{1}{3} \Lambda   \right) z^{-2/3} \\[6pt]
+ \left(  \frac{2}{9} (\Lambda - \Gamma) (e^{\eta_1 L} - 1)  +  \frac{1}{3} i p Le^{\eta_1 L}  \left( - \frac{2}{3} \Gamma  -  \frac{1}{3} \Lambda   \right)  \right) z^{-4/3} + O(z^{-5/3}). 
\end{multline}

On the other hand, we have 
$$
e^{\lambda_1 L + \tlambda_1 L + \lambda_3 L + \tlambda_3 L} = e^{-\lambda_2 L - \tlambda_2 L } \mathop{=}^{\eqref{lambda2-tlambda2}}  1 - \frac{1}{3} i p L  z^{-2/3}   +  O(z^{-4/3}). 
$$
and, by \Cref{th:lambda_asym},  
\begin{multline*}
\frac{1}{(\lambda_3 - \lambda_2)  (\tlambda_3 -
\tlambda_2)} = \frac{z^{-2/3}}{(\mu_3 - \mu_2) (\tmu_3 - \tmu_2)} - \frac{1}{3} \frac{z^{-4/3}}{(\mu_3 - \mu_2) (\tmu_3 - \tmu_2)} \Big( \frac{1}{\mu_2 \mu_3} + \frac{1}{\tmu_2 \tmu_3} \Big) + O (z^{-5/3}) \\[6pt]
= \frac{1}{3}z^{-2/3} - \frac{1}{9} z^{-4/3} + O (z^{-5/3}).
\end{multline*}
It follows that 
\be \label{sum2}
\frac{1}{(\lambda_3 - \lambda_2) (\tlambda_3 - \tlambda_2)} 
e^{\lambda_1 L + \tlambda_1 L + \lambda_3 L + \tlambda_3 L}
 =  \frac{1}{3}z^{-2/3}  \Big( 1 - \frac{1}{3} z^{-2/3} - \frac{1}{3} i p L z^{-2/3} \Big) + O (z^{-2}). 
\ee

Combining \eqref{sum1} and \eqref{sum2} yields \eqref{lem-B-cl1} for $z$ positive. 

We next deal with \eqref{lem-B-cl2}. 
When $e^{\eta_1 L} = 1$, we have 
\begin{multline}\label{sum1-1}
(T_1 + T_2)  + \Big(e^{\lambda_2 L + \tlambda_2 L + \eta_1 L} - 1\Big) (T_1 + \hat T_3) \mathop{=}^{\eqref{T1-T2}} \Big(e^{\lambda_2 L + \tlambda_2 L } - 1 \Big) (T_1 + \hat T_3)  + O(z^{-7/3}) \\[6pt]
\mathop{=}^{\eqref{lambda2-tlambda2}, \eqref{T1-T3}}  \left(\frac{1}{3} i p L  z^{-2/3}   - \frac{1}{9} i p L z^{-4/3} \right)  \left( \left( - \frac{2}{3} \Gamma  -  \frac{1}{3} \Lambda   \right) z^{-2/3}   
+ \frac{2}{9}(\Lambda -  \Gamma) z^{-4/3}  \right) + O(z^{-7/3}) \\[6pt]
= - \frac{(2 \Gamma + \Lambda)}{9} z^{-4/3} i p L \left(1 - \frac{1}{3} z^{-2/3} \right) \left(1 + \frac{2 (\Gamma - \Lambda)}{3 (2 \Gamma + \Lambda)} z^{-2/3} \right)+ O(z^{-7/3}). 
\end{multline}

Combining \eqref{sum2} and \eqref{sum1-1} yields \eqref{lem-B-cl2} for $z$ positive. 


\medskip
The conclusion in the case where $z$ is large and negative can be derived from the
case where $z$ is positive and large as follows.
Define, for $(z, x) \in  \mR \times (0, L)$, with large $|z|$,
\[
M(z, x) =  \frac{\sum_{j=1}^3(e^{\lambda_{j+1} L } - e^{\lambda_{j} L
})e^{\lambda_{j+2} x} }{\sum_{j=1}^3 (\lambda_{j+1} - \lambda_{j}) e^{-\lambda_{j+2} L
}}.
\]
Then
\[
B(z, x) = M(z, x) \overline{M (z-p, x)} \varphi_x(x).
\]
It is clear from the definition of $M$ that
\[
M(-z, x) = \overline{M(z, x)}.
\]
We then have
\[
B(-z, x) =M(-z, x) \overline{M (- z-p, x)} \varphi_x(x) =  \overline{ M(z, x)
\overline{M (z+ p, x)}  \, \overline{\varphi_x(x)}} .
\]
We thus obtain the result in the case where $z$ is negative and large  by taking the
conjugate of the corresponding expression for  large positive $z$ in which $\eta_j$
and $p$ are replaced by $-\eta_j$ and $-p$. The conclusion follows.
\end{proof}

\begin{remark} \rm The asymptotic expansion of the RHS of \eqref{lem-B-cl1} up to the order $|z|^{-4/3}$ was given in \cite[Lemma 3.4]{CKN-24} for $z \in \mR$ with large $|z|$, i.e., 
\be
 \int_0^L B(z, x) \diff x  = E |z|^{-4/3}   
+ O(|z|^{-5/3})
\mbox{ for $z\in \mR$ with large $|z|$}. 
\ee
The proof given here is more involved but in the same spirit.  In particular, we also rearrange various term differently from \cite{CKN-24} to better capture the cancellations. 
Assertion \eqref{lem-B-cl1} is sufficient to establish \Cref{thm-N-1} and \Cref{thm-N-2}. Nevertheless, both assertions \eqref{lem-B-cl1} and \eqref{lem-B-cl2} are needed in the proof of \Cref{thm-time} on the controllable time. 
\end{remark}

We next deal with $E$ and $F$. 

\begin{lemma} \label{lem-G} Let $k, \, l \in \N_*$, and $L>0$, and let $\eta_j$ for $j=1, 2, 3$ be defined by 
\begin{equation}
\eta_1 = - \frac{2 \pi i}{3 L} (2k + l), \quad \eta_2 = \eta_1 + \frac{2\pi i}{L} k, \quad
\eta_3 = \eta_2 + \frac{2\pi i}{L} l,
\end{equation}
Set 
\begin{equation}\label{lem-G-def-p}
p = \frac{(2k + l)(k-l)(2 l + k)}{3 \sqrt{3}(k^2 + kl + l^2)^{3/2}}.
\end{equation}
and assume that \eqref{def-L} holds. 
We have
\be
\Gamma = \Lambda   = - \frac{8 \pi^3}{L^3} i  k l (k+l). 
\ee
Consequently, 
\be
E =  - \frac{ \Gamma }{3} (e^{\eta_1 L} - 1) 
\quad \mbox{ and } \quad 
F = - \frac{\Gamma}{9}   i p Le^{\eta_1 L}.
\ee
\end{lemma}

\begin{proof}With $\gamma_j = L \eta_j/ (2 \pi i)$, we have
\[
\gamma_1 = - \frac{2k + l}{3}, \quad \gamma_2 = \frac{k-l}{3}, \quad \gamma_3 = \frac{k+ 2
l}{3}.
\]
It follows that
\begin{align*}
\frac{L^3}{(2 \pi i)^3}\sum_{j=1}^3 \eta_{j+2}^2(\eta_{j+1}-\eta_{j})  =  &   \sum_{j=1}^3
\gamma_{j+2}^2(\gamma_{j+1}-\gamma_{j})
=   \gamma_3^2 k + \gamma_1^2 l  - \gamma_2^2 (k+l)
\\[6pt] = &(\gamma_3^2  - \gamma_2^2) k  - (\gamma_2^2 - \gamma_1^2)l =  (k+l) k l,
\end{align*}
which yields
\begin{equation*}
\sum_{j=1}^3 \eta_{j+2}^2(\eta_{j+1}-\eta_{j})  = - \frac{8 \pi^3 i  k l (k+l)}{L^3}.
\end{equation*}

We also have
\begin{align} 
\sum_{j=1}^3  \frac{\eta_{j+1} - \eta_j}{\eta_{j+2}} =  &   \sum_{j=1}^3
\frac{\gamma_{j+1} - \gamma_j}{\gamma_{j+2}}
=  \frac{3k}{k + 2 l } - \frac{3 l }{2 k + l} - \frac{3 (k+l)}{k -l} = - \frac{2 7 k l
(k+l)}{(k+2l) (2 k + l) (k-l)}.
\end{align}
From \eqref{lem-G-def-p} and \eqref{def-L}, we have
\[
 \frac{p}{(k -l) (k+ 2l) (2l + k)} = \Big(\frac{2 \pi}{ 3 L} \Big)^3.
\]
We derive that
\[
i p \sum_{j=1}^3  \frac{\eta_{j+1} - \eta_j}{\eta_{j+2}} 
 = - \frac{8 \pi^3}{L^3}  i  kl (k+l).
\]
The proof is complete.
\end{proof}

As a consequence of  \Cref{lem-1,lem-B}, we obtain

\begin{lemma}\label{lem-dir} Let $k, \, l \in \N_*$, and $L>0$, and let $\eta_j$ be defined by \eqref{eta-kdv} for $j=1, 2, 3$.  
Set 
\begin{equation}
p = \frac{(2k + l)(k-l)(2 l + k)}{3 \sqrt{3}(k^2 + kl + l^2)^{3/2}}.
\end{equation}
Assume  that \eqref{def-L} holds. Define $\varphi$ by \eqref{def-varphi}. Let $u \in L^2(0, + \infty)$ and let $y \in C([0, + \infty); L^2(0, L))
\cap L^2_{\loc}\big([0, +\infty); H^1(0, L) \big)$ be the unique solution of
\eqref{sys-y} and \eqref{IC-y}.  Assume that $u(t) = 0$  and  $y(t, \cdot) = 0$ for large $t$.  We have
\begin{multline*}
\int_{0}^{+\infty} \int_0^L  |y(t, x)|^2 \varphi_x(x) e^{- i p t} \diff x \diff t  \\[6pt]
= \left\{ \begin{array}{c}
\dsp \int_{\mR}  \hu(z) \overline{\hu(z - p )} \Big(\frac{F}{|z|^{2}}  + \frac{F_1}{|z|^{8/3}}+ O(|z|^{-9/3})
\Big) \diff z \mbox{ if } e^{\eta_1 L } = 1, \\[6pt]
\dsp \int_{\mR}  \hu(z) \overline{\hu(z - p )} \Big(\frac{E}{|z|^{4/3}}  + \frac{E_1}{|z|^{2}}  + O(|z|^{-7/3})
\Big) \diff z \mbox{ if } e^{\eta_1 L } \neq 1. 
\end{array} \right. 
\end{multline*}
\end{lemma}

Here is an important consequence of \Cref{lem-dir}.

\begin{proposition} \label{pro-monotone-1} Let $k, \, l \in \N_*$, and $L>0$, and let $\eta_j$ be defined by \eqref{eta-kdv} for $j=1, 2, 3$.  
Set 
\begin{equation}
p = \frac{(2k + l)(k-l)(2 l + k)}{3 \sqrt{3}(k^2 + kl + l^2)^{3/2}}.
\end{equation}
Assume  that \eqref{def-L} holds. Define $\varphi$ by \eqref{def-varphi}.
Let $T>0$, and let $u \in L^2(0, + \infty)$ and let $y \in C([0, + \infty); L^2(0, L))
\cap L^2_{\loc}\big([0, +\infty); H^1(0, L) \big)$ be the unique solution of
\eqref{sys-y} and \eqref{IC-y}.  Assume that $u(t) = 0$ for $t > T$,  and  $y(t, \cdot) = 0$ for large $t$. 
Then
\be\label{pro-monotone-dir}
\int_{0}^{\infty} \int_0^L |y (t, x)|^2 e^{-ipt }\varphi_x(x) \diff x \diff t 
= \left\{
\begin{array}{c}
 \| u \|_{H^{-1}(\mR)}^2 \big(F + O(T^{1/3}) \big) 
 \mbox{ if } e^{\eta_1 L} = 1  \\[6pt] 
\| u \|_{H^{-2/3}(\mR)}^2 (E + O(T^{1/3}))\big)  \mbox{ if } e^{\eta_1 L } \neq 1. 
\end{array} \right. 
\ee
\end{proposition}

In \Cref{pro-monotone-1}, $\| u \|_{H^{-s}(\mR)}$ denotes the $H^{-s}(\mR)$-norm of the extension of $u$ by $0$ on $(-\infty, 0]$.

\begin{proof} By \Cref{pro-Gen},
\[
\hu G/ H \mbox{ is an entire function}.
\]
By \cite[Lemma 2.6]{CKN-24},  $G$ and $H$ are entire functions. One can show that the number of  common roots of $G$ and $H$ in $\mC$
is finite, see \cite[Lemma B.2]{CKN-24}. Let $z_1, \dots, z_k$ be the distinct common roots of $G$ and $H$ in $\mC$. There exist $m_1, \dots, m_k \in \N_*$ such that, with
\[
\Gamma(z) = \prod_{j=1}^k (z- z_j)^{m_j} \quad  \mbox{ in } \mC,
\]
the following two functions are entire
\begin{equation}
\quad \cG(z) : = \frac{G(z)}{\Gamma(z)} \quad \mbox{ and } \quad \cH(z) : = \frac{H(z)}{\Gamma(z)},
\end{equation}
and $\cG$ and $\cH$ have no common roots. Since
\[
 \hu \cG / \cH = \hu G/ H
\]
which is an  entire function, it follows that the function $v$ defined by
\begin{equation}
v(z) = \hu(z)/ \cH(z) = \hu(z) \frac{\Gamma(z)\Xi(z)}{\det Q(z)} \mbox{ in } \mC
\end{equation}
is also an entire function. We also have, see \cite[(3.41)]{CKN-24}, 
\begin{equation}\label{pro-monotone-1-vvv-p3}
|v(z)| \le C_\eps e^{(T+   \eps) |z|} \mbox{ for } z \in \mC.
\end{equation}

Since $\cH$ is an entire function, there exists $\gamma > 0$ such that
\begin{equation}\label{pro-monotone-H'}
\cH'''(z + i \gamma) \neq 0 \mbox{ for all } z \in \mR.
\end{equation}
Fix such an $\gamma$ and denote 
\be
\cH_\gamma (z) = \cH(z + i \gamma) \mbox{ for } z \in \mC.
\ee

We next derive some asymptotic behaviors for $\cH_\gamma$ and its derivatives. Since $\sum_{j=1}^3 \lambda_j = 0$, it follows from \eqref{eq-defQ} that 
\[
 \det Q =
(\lambda_2-\lambda_1)e^{ - \lambda_3 L} +  (\lambda_3-\lambda_2)e^{ - \lambda_1 L} +   (\lambda_1-\lambda_3)e^{ - \lambda_2 L}.
\]
We use the convention in  \Cref{th:lambda_asym}.  Thus, by \Cref{th:lambda_asym}, for a fixed $\beta \in \mR$, and large positive $z$, 
\begin{equation}\label{eq-H-asym}
\cH(z + i \beta) = \frac{\det Q(z +  i \beta)}{\Xi(z + i \beta)\Gamma(z + i \beta)} = \kappa 
z^{-2/3- \sum_{i=1}^k m_j}e^{ - \mu_1 L z^{1/3}}\big(1+O(z^{-1/3})\big), 
\end{equation}
where 
\[
\kappa = -  \frac{1}{ (\mu_2 - \mu_1) (\mu_1 - \mu_3)}. 
\]
We can also compute the asymptotic expansion of $\cH'(z + i \beta) $, either  by explicitly computing the asymptotic behavior of $\lambda_j'(z + i \beta)$ for large positive $z$ (formally, one just needs to take the derivative of \eqref{eq-H-asym} with respect to $z$). Similarly, we can also compute the asymptotic expansion of $\cH''(z + i \beta)$ and $\cH'''(z + i \beta)$.

We next only deal with the case $e^{\eta_1 L} = 1$. The proof in the other case is in the same spirit (see also \cite{CKN-24}).  
We get, for large positive $z$, 
\[
 \cH'''(z + i \beta) = - \frac{\mu_1^3 L^3 }{27} z^{-2}  \kappa 
z^{-2/3- \sum_{i=1}^k m_j}e^{ - \mu_1 L z^{1/3}}\big(1+O(z^{-1/3})\big). 
\]
We then get 
\[
\lim_{z \in \mR, z \to + \infty} z (z + i ) \frac{\cH_\gamma'''(z)}{\cH(z)}  = \frac{1}{\sigma}: = - \frac{\mu_1^3 L^3}{27} =  \frac{i L^3}{27}.
\]

Similarly, we obtain
\[
\lim_{z \in \mR, z \to - \infty} z (z + i ) \frac{\cH_\gamma'''(z)}{\cH(z)}    = \frac{1}{\bar \sigma}.
\]
Moreover, we have
\begin{equation}\label{pro-monotone-H}
\big |\cH(z) z^{-1} -   \sigma (z + i)  \cH_\gamma'''(z) \big|  \le  C |(z+i)\cH_\gamma'''(z)| |z|^{-1/3} \mbox{ for large positive $z$},
\end{equation}
and
\begin{equation}\label{pro-monotone-H-*}
\big |\cH(z) z^{-1} -  \bar \sigma  (z + i) \cH_\gamma'''(z) \big| \le C |(z+i)\cH_\gamma'''(z)| |z|^{-1/3} \mbox{ for negative $z$ with large $|z|$}.
\end{equation}

Set
\begin{equation}\label{eq-def-w}
\hw(z) =  (z + i ) v(z) \cH_\gamma'''(z) = \hu(z) (z + i) \frac{\cH'''_\gamma(z)}{ \cH(z)}.
\end{equation}
Then 
\begin{multline}
\mbox{$\hw$ is an entire function and satisfies Paley-Wiener's conditions} \\[6pt]  \mbox{for the
interval $(-T - \eps, T + \eps)$ for all $\eps >0$}. 
\end{multline}

Then $\hw$ is an entire function and satisfies Paley-Wiener's conditions for the
interval $(-T - \eps, T + \eps)$ for all $\eps >0$. Indeed, this follows from the
facts  $|\hw(z)| \le C_\eps  |v(z)| e^{\eps |z|}$ for $z \in \mC$ by \Cref{th:lambda_asym}, $|v(z)| \le C_\eps e^{(T + \eps)
|z|}$ for $z \in \mC$ by \eqref{pro-monotone-1-vvv-p3}, $|\cH_\gamma '(z) v(z)| =
|\cH_\gamma'(z)\cH(z)^{-1} \hu(z)| \le |\hat u(z)|$ for  real $z$ with large $|z|$,
so that $\int_{\mR} |\hw|^2 < + \infty$.

We have, by \Cref{lem-B} and \cite[Lemma B.1]{CKN-24},  
\begin{equation}\label{pro-monotone-BB}
\left| \int_0^L B(z, x) \diff x \right| \le \frac{C}{(|z| + 1)^{2}} \mbox{ for } z \in
\mR.
\end{equation}

From \eqref{pro-monotone-H'}, \eqref{pro-monotone-H},  \eqref{pro-monotone-H-*},  and \eqref{pro-monotone-BB},
we derive from \Cref{lem-B} that
\begin{equation*}
\left|\hu (z) \overline{\hu (z - p)} \int_0^L B(z, x) \diff x - F |\sigma|^2  \hw (z) \overline{\hw (z - p)} \right|  \le  C \frac{|\hw(z)|
|\hw (z-p)|}{(1 + |z|)^{1/3}} \mbox{ for } z \in \mR, 
\end{equation*}
which yields 
\begin{multline} \label{pro-monotone-B}
\left|\int_{\mR} \hu (z) \overline{\hu (z - p)} \int_0^L B(z, x) \diff x \diff z  - F
|\sigma|^{2}
\int_{\mR} \hw (z) \overline{ \hw(z-p)} \diff z  \right| \\[6pt]
\le C \int_{\mR} \frac{|\hat w(z)| |\hat w (z-p)|}{(1 + |z|)^{1/3}} \, dz 
\le C \int_{\mR} \frac{|\hat w(z)|^2}{(1 + |z|)^{1/3}} \, dz. 
\end{multline}

By \Cref{lem-interpolation} below, we have 
\be \label{pro-monotone-1-claim1}
\int_{\mR} \frac{|\hat w(z)|^2}{(1 + |z|)^{1/3}} \, dz  \le C T^{1/3} \| w \|_{L^2(\mR)}^2. 
\ee
Using the fact
\[
\int_{\mR} \hw (z) \overline{ \hw(z-p) } \diff z  =  \int_{\mR} |w(t)|^2 e^{-itp} \diff t
=
 \int_{-T}^{T} |w(t)|^2 e^{-itp} \diff t,
\]
we obtain from \eqref{pro-monotone-B} and \eqref{pro-monotone-1-claim1} that 
\[
\int_{\mR} \hu (z) \overline{\hu (z - p)} \int_0^L B(z, x) \diff x \diff z = 
|\sigma|^{2}
\int_{-T}^T |w(t)|^2 \Big(F +  O(1) T^{1/3} \Big) \diff  t.
\]
The conclusion follows by noting that
$$
\int_{\mR} \left|  |\sigma|^2 |\hw (z)|^2  -  \frac{|\hu (z)|^2}{(1 + |z|^2)} \right| \, dz \mathop{\le}^{\eqref{pro-monotone-H}, \eqref{pro-monotone-H-*}} \int_{\mR} \frac{|\hw (z)|^2}{(1 + |z|)^{1/3}} \mathop{\le}^{\Cref{lem-interpolation}} C T^{1/3} \int_{-T}^T |w(t)|^2  \, d t.
$$
The proof is complete. 
\end{proof}

In the proof of \Cref{pro-monotone-1}, we used the following lemma.

\begin{lemma} \label{lem-interpolation} Let $0< \alpha < 1$, $T> 0$,  and let $w \in L^2(\mR)$ be such that $\supp w \subset [-T, T]$. There exists a positive constant $C$ depending only on $\alpha$ such that 
\be 
\int_{\mR} \frac{|\hat w(z)|^2}{(1 + |z|)^{\alpha}} \, dz  \le C T^{\alpha} \| w \|_{L^2(\mR)}^2. 
\ee
\end{lemma}

\begin{proof} We have, for $m >0$,  
$$
 \int_{\mR} \frac{|\hat w(z)|^2}{(1 + |z|)^{\alpha}} \, dz 
\le  C \int_{|z| \le m} \frac{1}{(1+|z|)^{\alpha}}|\hw (z)|^2  \diff  z + C m^{-\alpha}
\int_{|z| > m} |\hw (z)|^2 \diff z, 
$$
for some positive constant $C$ depending only on $\alpha$. Since, for $z \in \mR$,
\[
|\hw(z)| \le C\| w\|_{L^1(\mR)} = C \|w \|_{L^1(-T, T)} \le C T^{1/2} \|w \|_{L^2(\mR)},
\]
we derive that
\begin{multline*}
\int_{\mR} \frac{|\hat w(z)|^2}{(1 + |z|)^{1/3}} \, dz \le C \left( T \int_{|z| < m} \frac{1}{(1 + |z|)^{\alpha}} \, dz  + m^{- \alpha} \right) \| w\|_{L^2(\mR)}^2\\[6pt]
\le C \Big( T m^{1 - \alpha} + m^{- \alpha} \Big) \| w\|_{L^2(\mR)}^2.
\end{multline*}
The conclusion follows by choosing $m = 1/T$. 
\end{proof}

We next present a
consequence of  \Cref{pro-monotone-1}. 
Multiplying \eqref{pro-monotone-dir} by  $\overline{E}$ if $ e^{\eta_1 L} \neq 1$ and by $\overline{F}$ if $ e^{\eta_1 L} = 1$, we obtain the following result from \Cref{pro-monotone-1}. 

\begin{corollary} \label{corollary-dir-1} Assume that the assumptions of \Cref{pro-monotone-1}
hold.  For any  (real) $u \in
L^2(0, + \infty)$ with $u(t) = 0$ for $t > T$ and $y(t,
\cdot) = 0$ for large $t$ where $y$ is the unique solution of \eqref{sys-y} and \eqref{IC-y}, 
we have
\begin{equation}\label{pro-monotone-CSQ1-*}
\int_{0}^\infty \int_0^{+\infty} y^2(t, x) \Psi_x(t, x) \diff x \diff t  = \left\{\begin{array}{c}  \| u
\|_{H^{-1}(\mR)}^2 (|F|^2 + O(T^{1/3}))  \mbox{ if } e^{\eta_1 L} = 1, \\[6pt]
\| u
\|_{H^{-2/3}(\mR)}^2 (|E|^2 + O(T^{1/3}))  \mbox{ if } e^{\eta_1 L} \neq 1, 
\end{array} \right. 
\end{equation}
where $\Psi$ is given by 
\begin{equation}\label{def-Psi-1}
\Psi(t, x) = \left\{ \begin{array}{c} \Re \Big( \bar F \varphi_x(x)e^{-i p t} \Big) \mbox{ if } e^{\eta_1 L} = 1, \\[6pt]
\Re  \Big( \bar E \varphi_x(x) e^{-i p t} \Big) \mbox{ if } e^{\eta_1 L} \neq 1. 
\end{array} \right. 
\end{equation}
Consequently, there exists $T_*> 0$ such that if in addition that $T < T_*$ then 
\begin{equation}\label{pro-monotone-CSQ1}
\int_{0}^\infty \int_0^{+\infty} y^2(t, x) \Psi_x(t, x) \diff x \diff t \ge  \left\{\begin{array}{c}  \frac{1}{2} |F|^2 \| u
\|_{H^{-1}(\mR)}^2 \mbox{ if } e^{\eta_1 L} = 1, \\[6pt] \frac{1}{2} |E|^2 \| u
\|_{H^{-2/3}(\mR)}^2 \mbox{ if } e^{\eta_1 L} \neq 1. 
\end{array} \right. 
\end{equation}
\end{corollary}

\begin{remark} \rm \Cref{corollary-dir-1} in the case $e^{\eta_1 L } \neq 1$ was already established in \cite{CKN-24} in a very closely related form (see \cite[Corollary 3.7]{CKN-24}). 
\end{remark}

\section{Proof of \Cref{thm-N-1,thm-N-2}} \label{sect-thm-N}

We first recall the following two results from \cite{Ng-KdV-D} which will be used in the proof of \Cref{thm-N-1,thm-N-2}. Here is the first one, which is a special case of  \cite[Proposition 3.1]{Ng-KdV-D}. 

\begin{lemma}\label{lem-kdv1} Let $L>0$,  $T>0$, $h \in L^2 (0, T)$,  $ f \in L^1\big((0, T); L^2(0, L)\big)$, and $y_0 \in L^2(0, L)$. There exists a unique solution $y \in X_T$  of the system
\begin{equation}\label{sys-y-LKdV}\left\{
\begin{array}{cl}
y_t + y_x  + y_{xxx}  = f &  \mbox{ in } (0, T) \times  (0, L),
\\[6pt]
y(\cdot, 0) = 0,  \;  y(\cdot, L) = 0, \;  y_x(\cdot, L)   = h & \mbox{ in } (0,
T), \\[6pt]
y(0, \cdot)  = y_0 & \mbox{ in } (0, L).
\end{array}\right.
\end{equation}
Moreover, for $x \in [0, L]$,
\begin{multline}\label{pro-kdv1-cl1} 
\| y\|_{X_T} + \| y(\cdot, x)\|_{H^{1/3}(0, T)} + \| y_x(\cdot, x)\|_{L^2(0, T)}  \\[6pt]
\le C_{T, L} \Big( \| y_0\|_{L^2(0, L)} + \| f\|_{L^1\big( (0, T);  L^2(0, L) \big)} + \|
h \|_{L^2(0, T)}\Big),  
\end{multline}
and
\begin{multline} \label{pro-kdv1-cl2} 
\| y(\cdot, x)\|_{L^2(0, T)} + \| y_x(\cdot, x)\|_{[H^{1/3}(0, T)]^*} \\[6pt]
\le C_{T, L} \Big( \| y_0\|_{[H^{1}(0, L)]^*} + \| f\|_{L^1\big( (0, T);  [H^{1}(0, L)]^* \big)} + \|
h \|_{[H^{1/3}(0, T)]^*}\Big),
\end{multline}
where  $C_{T, L}$ denotes a positive constant  independent of $x$,  $y_0$, $f$, and $h$. 
\end{lemma}

Here is the second one. 

\begin{lemma}[\mbox{\cite[Lemma 8.1]{Ng-KdV-D}}] \label{lem-kdv2} Let $L > 0$ and $T > 0$, and let  $f \in L^2 \big( (0, T); H^1(0, L) \big)$ and $g \in L^2\big((0, T) \times (0, L) \big)$. Then 
\be\label{lem-interpolationL1L2-st1}
\| f g \|_{L^1\big((0, T); L^2(0, L)\big)} \le C \| f\|_{L^2 \big((0, T) \times (0, L) \big)}^{\frac{1}{2}}  \| f\|_{L^2 \big((0, T); H^1 (0, L) \big)}^{\frac{1}{2}}  \| g\|_{L^2 \big((0, T) \times (0, L) \big)} 
\ee
and 
\begin{multline}\label{lem-interpolationL1L2-st11}
\| (fg)_x\|_{L^1((0, T); [H^1(0, L)]^*)} \\[6pt]
\le C \Big(\| f g \|_{L^1((0, T); L^2(0, T))} + \| (f g)(\cdot, 0) \|_{L^1(0, T)} + \| (f g)(\cdot, L) \|_{L^1(0, T)} \Big),  
\end{multline}
for some positive constant $C$ depending only on $L$. Consequently, for $f, g \in L^2 \big( (0, T); H^1(0, L) \big)$, 
\be\label{lem-interpolationL1L2-st2}
\| f g_x \|_{L^1\big((0, T); L^2(0, L)\big)} \le C \| f\|_{L^2 \big((0, T) \times (0, L) \big)}^{\frac{1}{2}}  \| f\|_{L^2 \big((0, T); H^1 (0, L) \big)}^{\frac{1}{2}}  \| g_x\|_{L^2 \big((0, T) \times (0, L) \big)}.
\ee
\end{lemma}

We are ready to give the proof of \Cref{thm-N-1,thm-N-2}. 

\begin{proof}[Proofs of \Cref{thm-N-1,thm-N-2}] Let $\eps_0$ be a small positive constant, which depends only on $L$,   and is determined
later.  We prove \Cref{thm-N-1,thm-N-2} with $\Phi$ given by $\Psi(0, \cdot)$ where $\Psi$ is the function defined by \eqref{def-Psi-1}. Assume that there exists a solution $y
\in C\big([0, + \infty); L^2(0, L) \big) \cap L^2_{\loc}\big([0, + \infty); H^1(0, L)
\big) $ of
\eqref{sys-KdV-N}  with $y(t, \cdot) =0$ for $t \ge T/2$,  for some $0< \eps < \eps_0$, for some $u \in L^2(0, +
\infty)$  with   $\supp u \subset [0, T/2]$ and $\| u
\|_{H^{s}(0, T)} < \eps_0$ where $s = 1/2$ if $2 k + l \not \in 3 \N_*$ and $s = 7/6$ if $2 k + l \not \in 3 \N_*$. We also assume that $T \le 1$. 

Using the fact $y(t, \cdot) = 0$ for $t \ge T/2$, we have, for $\eps_0$ small, 
\begin{equation}\label{thm-NL-yyy}
\| y\|_{L^2 \big(\mR_+; H^1(0, L) \big)} \le C \Big( \| y_0 \|_{L^2(0, L)} + \| u
\|_{L^2(\mR_+)} \Big), 
\end{equation}
which in turn implies,  by \Cref{lem-kdv1,lem-kdv2}, 
\begin{equation}\label{thm-NL-yyy-1}
\| y\|_{L^2 \big( \mR_+ \times (0, L) \big)} \le C \Big( \| y_0 \|_{L^2(0, L)} + \| u
\|_{[H^{1/3}(\mR_+)]^*} \Big).
\end{equation}
Here and in what follows, $C$ denotes a positive constant depending only on $L$ and $T$. 

Let $y_1$ and  $y_2$ be the solutions of the following systems
\begin{equation} \label{thm-NL-y1}\left\{
\begin{array}{cl}
y_{1, t}  + y_{1, x}  + y_{1, xxx}  = 0 &  \mbox{ in } 
\mR_+ \times (0, L), \\[6pt]
y_1( \cdot, 0) =   y_{1}(\cdot, L) = 0 & \mbox{ in }  \mR_+,\\[6pt]
y_{1, x}(\cdot, L) = u& \mbox{ in }  \mR_+,\\[6pt]
y_{1}(0 , \cdot) = 0  & \mbox{ in }  (0, L), 
\end{array}\right.
\end{equation}
\begin{equation}\label{thm-NL-y2}\left\{
\begin{array}{cl}
y_{2, t}  + y_{2, x}  + y_{2, xxx}  + y_{1} y_{1, x}= 0 &  \mbox{ in } 
\mR_+ \times (0, L), \\[6pt]
y_2( \cdot, 0) =   y_{2, x}(\cdot, L) = y_{2}(\cdot, L)= 0 & \mbox{ in }  \mR_+,\\[6pt]
y_{2}(0 , \cdot) = 0 & \mbox{ in }  (0, L). 
\end{array}\right.
\end{equation}

Then, by \Cref{lem-kdv1},   
\be \label{thm-NL-y1-1}
\|y_1\|_{X_T} \le C \| u\|_{L^2(\mR_+)} 
\ee
and
\be\label{thm-NL-y1-2}
\| y_1\|_{L^2 \big((0, T) \times (0, L) \big)} \le C    \| u\|_{[H^{1/3}(\mR_+)]^*},  
\ee
which in turn imply, by \Cref{lem-kdv1,lem-kdv2}, 
\be \label{thm-NL-y2-1}
\|y_2\|_{X_T} \le C \| y_{1} y_{1, x} \|_{L^1((0, T); L^2(0, L))} \le C \| u\|_{[H^{1/3}(\mR_+)]^*}^{1/2} \| u\|_{L^2(\mR_+)}^{3/2}, 
\ee
and
\be\label{thm-NL-y2-2}
\| y_2\|_{L^2 \big((0, T) \times (0, L) \big)} \le C \| y_1^2 \|_{L^1((0, T); L^2(0, L))}   \le C   \| u\|_{[H^{1/3}(\mR_+)]^*}^{3/2} \| u\|_{L^2(\mR_+)}^{1/2}. 
\ee

Set 
$$
\dy = y - y_1 - y_2 \mbox{ in } \mR_+ \times (0, L).  
$$
We have 
\begin{equation*}\left\{
\begin{array}{cl}
\dy_t  + \dy_x  + \dy_{xxx}  + y \dy_x + \dy (y_1 + y_2)_x  + (y_1 y_2)_x + y_2 y_{2, x}  = 0&  \mbox{ in } 
\mR_+ \times (0, L), \\[6pt]
\dy(\cdot, 0) =  \dy_{x} (\cdot, L) = \dy(\cdot, L)  = 0 & \mbox{ in } \mR_+, \\[6pt]
\dy(0, \cdot) = \eps \Phi & \mbox{ in } (0, L).
\end{array}\right.
\end{equation*}
Applying   \Cref{lem-kdv1,lem-kdv2},  we derive  that
\begin{multline}\label{thm-NL-dy-1}
\|\dy\|_{X_T} \le C \Big( \| y_{1} y_{2, x} \|_{L^1((0, T); L^2(0, L))} + \| y_{2} y_{1, x} \|_{L^1((0, T); L^2(0, L))}  + \| y_{2} y_{2, x} \|_{L^1((0, T); L^2(0, L))}  \Big)\\[6pt]
  \mathop{\le}^{\eqref{thm-NL-y1-1}-\eqref{thm-NL-y2-2}} C  \| u\|_{[H^{1/3}(\mR_+)]^*} \| u\|_{L^2(\mR_+)}^{2} + C \eps,
\end{multline}
and 
\begin{multline}\label{thm-NL-dy-2}
\|\dy\|_{L^2\big( (0, T) \times (0, L) \big)} \le C \Big( \|y \dy_x \|_{L^1((0, T); L^2(0, L))} + \|(y_1 + y_2) \dy_x \|_{L^1((0, T); L^2(0, L))}  \Big) \\[6pt]
+ C \Big(\|(y_1 + y_2) \dy \|_{L^1((0, T); L^2(0, L))} +  \| y_1 y_2 \|_{L^1((0, T); L^2(0, L))} +  \|y_2^2\|_{L^1((0, T); L^2(0, L))} \Big) \\[6pt]\mathop{\le}^{\eqref{thm-NL-yyy}-\eqref{thm-NL-yyy-1}, \eqref{thm-NL-y1-1}-\eqref{thm-NL-dy-1}}  C  \| u\|_{[H^{1/3}(\mR_+)]^*}^{3/2} \| u\|_{L^2(\mR_+)}^{5/2} + C \| u\|_{[H^{1/3}(\mR_+)]^*}^{2} \| u\|_{L^2(\mR_+)} 
+ C \eps 
 \\[6pt]
\le C    \| u\|_{[H^{1/3}(\mR_+)]^*}^{3/2} \| u\|_{L^2(\mR_+)}^{3/2} + C \eps. 
\end{multline}
In \eqref{thm-NL-dy-1}, we absorbed the contribution in the RHS of 
$$
\| y \dy_x\|_{L^1((0, T); L^2(0, L))} + \|\dy (y_1 + y_2)_x \|_{L^1((0, T); L^2(0, L))},
$$ 
which is bounded by $C \eps_0 \| \dy \|_{X_T}$. Combining  \eqref{thm-NL-y2-2} and \eqref{thm-NL-dy-2} yields
\be\label{thm-NL-y1-0}
\|y - y_1 \|_{L^2\big((0, T) \times (0, L) \big)} \le C    \| u\|_{[H^{1/3}(\mR_+)]^*}^{3/2} \| u\|_{L^2(\mR_+)}^{1/2} + C \eps. 
\ee
Since $y =0$ for $t \ge T/2$ and $u = 0$ for $t \ge T/2$, after using the regularizing effect of the linear KdV equation, 
and considering the projection into $\cM_N^\perp$,  we derive that
\be\label{thm-NL-y1-***}
\|y_1 (T, \cdot) \|_{L^2(0, L)} \le C    \| u\|_{[H^{1/3}(\mR_+)]^*}^{3/2} \| u\|_{L^2(\mR_+)}^{1/2} + C \eps. 
\ee

Since $y_1(T, \cdot) \in \cM_{N}^\perp$, it follows that there exists  $u_1 \in L^2(0, T)$ such that 
\be\label{thm-NL-yyy1}
\| u_1 \|_{L^2(0, T)} \le C \| y_1(T, \cdot) \|_{L^2(0, L)}
\ee
and  the solution $\ty_1 \in X_{T} $ of the  system 
\begin{equation*}\left\{
\begin{array}{cl}
\ty_{1, t}  + \ty_{1, x}  + \ty_{1, xxx}   = 0 &
\mbox{ in } (0, T) \times (0, L), \\[6pt]
\ty_1(\cdot, 0) = \ty_{1}(\cdot, L)  = 0 & \mbox{ in }  (0, T), \\[6pt]
\ty_{1, x}(\cdot, L) =u_1 & \mbox{ in } (0, T), \\[6pt]
\ty_1(0 , \cdot) = 0 & \mbox{ in } (0, L),
\end{array}\right.
\end{equation*}
satisfies 
$$
\ty_1(T, \cdot) = - y_1(T, \cdot).
$$
Using  \eqref{thm-NL-y1-***}, we derive from \eqref{thm-NL-yyy1} that 
\be \label{thm-NL-u1}
\| u_1 \|_{L^2(0, T)} \le C   \| u\|_{[H^{1/3}(\mR_+)]^*}^{3/2} \| u\|_{L^2(\mR_+)}^{1/2} + C \eps, 
\ee
which in turn implies, by \Cref{lem-kdv1}, 
\be \label{thm-NL-ty1}
\| \ty_1 \|_{X_T} \le C  \| u\|_{[H^{1/3}(\mR_+)]^*}^{3/2} \| u\|_{L^2(\mR_+)}^{1/2} + C \eps. 
\ee

Set 
$$
\ty = y_1 + \ty_1 \mbox{ in } (0, T) \times (0, L). 
$$
Then 
\begin{equation}\label{thm-NL-ty}\left\{
\begin{array}{cl}
\ty_{t}  + \ty_{x}  + \ty_{xxx}   = 0 &
\mbox{ in } (0, T) \times (0, L), \\[6pt]
\ty(\cdot, 0) = \ty(\cdot, L)  = 0 & \mbox{ in }  (0, T), \\[6pt]
\ty_x (\cdot, L) = u + u_1 & \mbox{ in } (0, T), \\[6pt]
\ty(0 , \cdot) =  \ty(T, \cdot) =0  & \mbox{ in } (0, L). 
\end{array}\right.
\end{equation}

We have 
\begin{multline}\label{thm-NL-yhy}
\| y - \ty \|_{L^2\big( (0, T) \times (0, L)\big)} \le \| y - y_1 \|_{L^2\big( (0, T) \times (0, L)\big)} + \| \ty_1 \|_{L^2\big( (0, T) \times (0, L)\big)} \\[6pt]
\mathop{\le}^{\eqref{thm-NL-y1-0}, \eqref{thm-NL-ty1}} C   \| u\|_{[H^{1/3}(\mR_+)]^*}^{3/2} \| u\|_{L^2(\mR_+)}^{1/2} + C \eps.
\end{multline}

Multiplying the equation of $y$ with $\Psi(t, x)$, integrating by parts on $[0, L]$, we have
\begin{equation}\label{thm-NL-B}
\frac{d}{dt} \int_{0}^L y (t, x) \Psi(t, x)\diff x - \frac{1}{2} \int_0^L y^2 (t, x)
\Psi_x (t, x) \diff x = 0.
\end{equation}
Integrating \eqref{thm-NL-B} from 0 to $T$ and using the fact $y(T, \cdot) = 0$ yield
\begin{equation}\label{thm-NL-id}
\int_{0}^L y_0 (x)  \Psi(0, x) \diff x  + \frac{1}{2} \int_0^{T} \int_0^L y^2 (t, x)
\Psi_x(t, x) \diff x \diff t  = 0.
\end{equation}
We have
\begin{multline}\label{thm-NL-id1}
 \|  y^2 - \ty^2 \|_{L^1\big( (0, T) \times (0, L) \big)} \le  \|  y - \ty\|_{L^2\big( (0, T) \times (0, L) \big)}\|  y + \ty\|_{L^2\big( (0, T) \times (0, L) \big)}\\[6pt] \mathop{\le}^{\eqref{thm-NL-yyy-1}, \eqref{thm-NL-yhy}} C   \| u\|_{[H^{1/3}(\mR_+)]^*}^{5/2} \| u\|_{L^2(\mR_+)}^{1/2} + C \eps_0 \eps.
\end{multline}
Combining   \eqref{thm-NL-id} and \eqref{thm-NL-id1} yields 
\be\label{thm-NL-key} 
\int_{0}^L y_0 (x)  \Psi(0, x) \diff x   + \frac{1}{2} \int_0^{T} \int_0^L \ty^2 (t, x)
\Psi_x(t, x) \diff x \diff t  \le C   \| u\|_{[H^{1/3}(\mR_+)]^*}^{5/2} \| u\|_{L^2(\mR_+)}^{1/2} + C \eps_0 \eps. 
\ee

The rest of the proof is divided into two cases. 

\bigskip 

\noindent $\bullet$ {\it Case 1: $2k + l \not \in 3 \N_+$.} Applying  \Cref{corollary-dir-1} to $\ty$  after noting \eqref{thm-NL-ty}, we derive from \eqref{thm-NL-key} that, for $\eps_0$ sufficiently small,   
$$
\eps \|\phi \|_{L^2(0, L)}^2 + \frac{1}{2} \alpha \| u + u_1 \|_{[H^{2/3}(\mR_+)]^*}^2 \le C  \| u\|_{[H^{1/3}(\mR_+)]^*}^{5/2} \| u\|_{L^2(\mR_+)}^{1/2} + C \eps_0 \eps. 
$$
Using \eqref{thm-NL-u1}, it follows that, for  sufficiently small $\eps_0$,  
$$
\alpha \| u  \|_{[H^{2/3}(\mR_+)]^*}^2 \le C \| u\|_{[H^{1/3}(\mR_+)]^*}^{5/2} \| u\|_{L^2(\mR_+)}^{1/2}.
$$
We derive that   
$$
\alpha \| u  \|_{[H^{2/3}(0, T/2)]^*}^2 \le C\| u
\|_{[H^{1/3}(0, T/2)]^*}^{5/2}  \| u
\|_{L^2(0, T/2)}^{1/2}.
$$

We have, by \Cref{lem-interpolationA} in the appendix, 
$$
\| u\|_{L^2(0, T/2)} \le C \| u\|_{[H^{2/3}(0, T/2)]^*}^{3/7} \| u\|_{H^{1/2}(0, T/2)}^{4/7}
$$
and 
$$
\| u\|_{[H^{1/3}(0, T/2)]^*} \le C \| u\|_{[H^{2/3}(0, T/2)]^*}^{5/7} \| u\|_{H^{1/2}(0, T/2)}^{2/7}.
$$
This yields 
$$
\| u\|_{[H^{1/3}(0, T/2)]^*}^{5/2}\| u\|_{L^2(0, T/2)}^{1/2} \le C  \| u\|_{[H^{2/3}(0, T/2)]^*}^2 \| u \|_{H^{1/2}(0, T/2)}.  
$$

So, for fixed sufficiently small $\eps_0$,
\[
u =0. 
\]
Hence $y(T, \cdot)  \not \equiv 0$. We have a contradiction.

\bigskip 
\noindent $\bullet$ {\it Case 2: $2k + l \in 3 \N_+$}

Applying \Cref{corollary-dir-1} to $\ty$  after noting \eqref{thm-NL-ty}, we derive from \eqref{thm-NL-key} that, for $\eps_0$ sufficiently small,   
$$
\eps \|\phi \|_{L^2(0, L)}^2 + \frac{1}{2} \alpha \| u + u_1 \|_{[H^{1}(\mR_+)]^*}^2 \le C  \| u\|_{[H^{1/3}(\mR_+)]^*}^{5/2} \| u\|_{L^2(\mR_+)}^{1/2} + C \eps_0 \eps. 
$$
Using \eqref{thm-NL-u1}, it follows that, for  sufficiently small $\eps_0$,  
$$
\alpha \| u  \|_{[H^{1}(\mR_+)]^*}^2 \le C \| u\|_{[H^{1/3}(\mR_+)]^*}^{5/2} \| u\|_{L^2(\mR_+)}^{1/2}.
$$
We derive that   
$$
\alpha \| u  \|_{[H^{1}(0, T/2)]^*}^2 \le C\| u
\|_{[H^{1/3}(0, T/2)]^*}^{5/2}  \| u
\|_{[L^2(0, T/2)]^*}^{1/2}.
$$

We have, by \Cref{lem-interpolationA} in the appendix, 
$$
\| u\|_{L^2(0, T/2)} \le C \| u\|_{[H^{1}(0, T/2)]^*}^{7/13} \| u\|_{H^{7/6}(0, T/2)}^{6/13}
$$
and 
$$
\| u\|_{[H^{1/3}(0, T/2)]^*} \le C \| u\|_{[H^{1}(0, T/2)]^*}^{9/13} \| u\|_{H^{7/6}(0, T/2)}^{4/13}.
$$
This yields 
$$
\| u\|_{[H^{1/3}(0, T/2)]^*}^{5/2}\| u\|_{L^2(0, T/2)}^{1/2} \le C  \| u\|_{[H^{1}(0, T/2)]^*}^2 \| u \|_{H^{7/6}(0, T/2)}.  
$$

So, for fixed sufficiently small $\eps_0$,
\[
u =0. 
\]
Hence $y(T, \cdot) \not \equiv 0$. We have a contradiction.

\medskip 
The proof is complete. 
\end{proof}

\section{The controllable time of the KdV system - Proof of \Cref{thm-time}}
\label{sect-thm-time}

\subsection{Preliminaries}

We begin this section with the following result. 

\begin{lemma}  \label{lem-FE} We have 
\be
\frac{E_1}{E} = \left\{ \begin{array}{cl} \dsp - \frac{1}{3} + \frac{\sqrt{3}}{18} p L - \frac{1}{6} i p L & \mbox{ if } e^{\eta_1 L } = e^{ 2 \pi i / 3}, \\[6pt]
\dsp - \frac{1}{3} - \frac{\sqrt{3}}{18} p L - \frac{1}{6} i p L & \mbox{ if } e^{\eta_1 L } = e^{ 4 \pi i / 3}.
\end{array} \right.
\ee
Consequently, 
$$
\Im (E_1/E) = - \frac{1}{6} p L \mbox{ if } e^{i \eta_1 L } \neq 1. 
$$
\end{lemma}

\begin{proof} We first consider the case  $e^{\eta_1 L } = e^{ 2 \pi i / 3}$. We have 
$$
\frac{e^{i \eta_1 L }}{e^{\eta_1 L} - 1} = \frac{1}{2} - i \frac{\sqrt{3}}{6}. 
$$
This implies 
$$
\frac{E_1}{E} = - \frac{1}{3} (1 + i p L)  + \frac{1}{3} i p L \left( \frac{1}{2} - i \frac{\sqrt{3}}{6}\right) = - \frac{1}{3} + \frac{\sqrt{3}}{18} p L - \frac{1}{6} i p L. 
$$

We next consider the case $e^{\eta_1 L } = e^{ 4 \pi i / 3}$. We have 
$$
\frac{e^{i \eta_1 L }}{e^{\eta_1 L} - 1} = \frac{1}{2} + i \frac{\sqrt{3}}{6}. 
$$
This implies 
$$
\frac{E_1}{E} = - \frac{1}{3} (1 + i p L) + \frac{1}{3} i p L \left( \frac{1}{2} + i \frac{\sqrt{3}}{6}\right) = - \frac{1}{3} - \frac{\sqrt{3}}{18} p L - \frac{1}{6} i p L. 
$$

The conclusion follows. 
\end{proof}

We next state a useful result concerning $\cM_N$.  

\begin{lemma} \label{lem-orthogonal} Let $\Psi$ be defined by \eqref{def-Psi} and let $c \in \mC \setminus \{0 \}$. Given $\tau \in \mR$, define 
$$
\Psi_1 (x) = \Re \big( c \Psi (\tau, x) \big) \quad \mbox{ and } \quad \Psi_2 = \Im \big(  c \Psi(\tau, x) \big) \mbox { for } x \in [0, L].
$$ 
Then 
$$
\int_{0}^L |\Psi_1 (x)|^2 \, dx = \int_0^L |\Psi_2 (x)|^2 \, dx \quad \mbox{ and } \quad \int_0^L \Psi_1 (x) \Psi_2 (x) \, dx = 0. 
$$
\end{lemma}

\begin{proof} For notational ease, we assume that $\tau = 0$. We have 
$$
c \Psi(t, x) =  \Big(\Psi_1 (x) + i \Psi_2(x) \Big) \Big(\cos (pt) - i \sin (pt) \Big).  
$$
This implies that 
\be
y(t, x): = \Re \big( c\Psi(t, x) \big) =  \cos (p t) \Psi_1 (x) + \sin (p t) \Psi_2 (x).  
\ee
Using \eqref{Psi-1} and \eqref{Psi-2}, we derive that  
\be \label{lem-orthogonal-p1}
\left\{\begin{array}{cl}
y_t(t, x) + y_{x}(t, x) + y_{xxx} (t, x) = 0 \mbox{ in } \mR \times [0, L], \\[6pt]
y(t, 0) = y(t, L) = y_x(t, 0) = y_x(t, L) = 0 \mbox{ in } \mR. 
\end{array} \right. 
\ee
Multiplying the equation of $y$ by $y$ and integrating by parts, we obtain 
\be \label{lem-orthogonal-p2}
\frac{d}{dt} \int_0^L |y(t, x)|^2 \, dx = 0 \mbox{ for } t \in \mR.  
\ee
The conclusion now follows from \eqref{lem-orthogonal-p1} and \eqref{lem-orthogonal-p2}. 
\end{proof}

\subsection{The key ingredient}

Here is the main ingredient of the proof of \Cref{thm-time}. 

\begin{proposition} \label{pro-time1} Let $L \in \cN_N$ be such that $L \not \in \{  L(2, 1), L(3, 1) \}$,  and let $k_m, l_m \in \N_*$ be such that $L =  L(k_m, l_m)$. There exists $\eps > 0$ such that  for every $\varphi \in M_m$, there exists $u_1 \in L^2(0, p_m/\pi - \eps)$ such that 
$$
y_1 (p_m/\pi - \eps) = 0 \quad \mbox{ and } \quad y_2 (p_m/\pi - \eps) = \varphi, 
$$
where $y_1$ and $y_2$ are the unique solutions of \eqref{eq:first_order} and \eqref{eq:second_order}, respectively, with zero initial data. 
\end{proposition}

\begin{proof}
Define, for $\nu >1$ and $0 < \beta < 1$, 
$$
\hat v(z)= \int_{0}^2 e^{-\frac{\nu}{1-(t-1)^2}}e^{-i \beta t z} \,d t. 
$$
Then 
$$
\hat v(z)= e^{- i \beta z} \int_{-1}^1 e^{-\frac{\nu}{1-t^2}}e^{-i \beta t z} \,d t = e^{- i \beta z} \hat v_1(z). 
$$
We have, by \cite[Lemma 4.3]{Tenenbaum-Tucsnak} \footnote{In \cite[Lemma 4.3]{Tenenbaum-Tucsnak}, $\beta \nu$ is fixed by a specific constant. Nevertheless, the proof there also gives the result stated here.}, 
\be
|\hat v(z)| = |\hat v_1(z)| \le  C e^{- \frac{\nu}{4}} e^{-(1+\delta) \sqrt{\beta \nu |z|}} \mbox{ for } z \in \mR. 
\ee

We have 
$$
 - \sqrt{\beta \nu |z|} + \frac{\sqrt{3}}{2} |z|^{1/3} = 0 \mbox{ if } (\beta \nu)^3 z = \frac{3}{4}, 
$$
and, with $\beta \nu^2 = (1, 617)^2 >  ( \frac{3}{4} + \frac{\sqrt{3}}{2})^2$,  
\be \label{pro-time1-01}
-\frac{\nu}{4} - \sqrt{\beta \nu |z|} + \frac{\sqrt{3}}{2} |z|^{1/3} \le (-1 - 2 \delta ) \nu \mbox{ if } z \ge  \nu^3 \mbox{ for some $\delta > 0$.}
\ee

We will choose 
\be \label{pro-time1-beta-nu}
\beta = T /2 \quad \mbox{ and } \quad \nu^2 =  \frac{(1.617)^2}{\beta} = 5.223 T^{-1}. 
\ee


We now consider the two cases $e^{i \eta_1 L } \neq 1$ and $e^{i \eta_1 L } = 1$ separately.

\medskip
\noindent $\bullet$ {\it Case 1}: $e^{i \eta_1 L } \neq 1$.  Fix $\gamma \neq 0$ such that 
\be \label{pro-time1-H1}
H'(z + i \gamma) \neq 0 \mbox{ for } z \in \mR. 
\ee
Denote $H_\gamma (z) = H(z + i \gamma)$ and  set 
\be \label{pro-time1-def-uw}
\hat u (z) =  \hat v(z) H(z) \quad \mbox{ and } \quad \hat w (z) = \frac{3}{\mu_3} \hat v(z) H_\gamma'(z). 
\ee
Recall that $H$ is defined in \eqref{def-GH}. One can check that $u \in L^2(\mR)$ with $\supp u \subset [0, T]$; moreover, concerning the linearized system, $u$ is a control which steers $0$ from time $0$ to $0$ at time $T$. This follows from the fact that  $
e^{i \beta \cdot} \hu$ and $e^{i \beta \cdot} \hu G/H$ are entire functions and satisfy the assumptions of Paley-Wiener's theorem concerning the support in $[-T/2, T/2]$ (see also \Cref{lem-form-sol} and \Cref{pro-Gen}).

By \Cref{th:lambda_asym}, we have
\be \label{pro-time1-recall}
H(z) = \frac{e^{\mu_3 |z|^{1/3} - \frac{1}{3 \mu_3 |z|^{1/3}}  }}{(\mu_3 - \mu_1) (\mu_2 - \mu_1) |z|^{2/3}(1 + O(|z|^{-2/3}))} \Big(1 + O(|z|^{-2/3}) \Big) \mbox{ for real $z$ with  large $|z|$}. 
\ee

We have, by \Cref{th:lambda_asym}, for real $z$ with large $|z|$, 
$$
\hat w(z) = \frac{\hat u (z)}{|z|^{2/3}} \left(1  - \frac{2}{\mu_3 |z|^{1/3}} + \frac{3}{\mu_3^2 |z|^{2/3}} + O(|z|^{-1}) \right).  
$$
This implies, by \eqref{pro-time1-H1}, for $z \in \mR$, 
$$
\frac{\hat u (z)}{|z|^{2/3}} =\hat w(z) \left(1 + \frac{2}{\mu_3 |z|^{1/3}} +  \frac{1}{\mu_3^2 |z|^{2/3}} + O (|z|^{-1}) \right). 
$$
We thus derive that 
\begin{multline*}
\hat u (z) \overline{\hat u (z-p)} \left(\frac{1}{|z|^{4/3}} + \frac{E_1}{E} \frac{1}{|z|^{2}}  + O (|z|^{-7/3}) \right)= \hat w (z) \overline{\hat w (z-p)} +  4 \Re \left(\frac{1}{\mu_3} \right) \frac{\hat w (z) \overline{\hat w (z-p)}}{|z|^{1/3}}  \\[6pt]
+  \left(\frac{1}{\mu_3 \bar \mu_3} + \frac{1}{\mu_3^2} + \frac{1}{\bar \mu_3^2} \right)
\frac{\hat w (z) \overline{\hat w (z-p)} }{|z|^{2/3}}  
+ \frac{E_1}{E}  \frac{\hat w (z) \overline{\hat w (z-p)}}{|z|^{2/3}} + \hat w (z) \overline{\hat w (z-p)} O(|z|^{-1}). 
\end{multline*}
Since 
$$
4 \Re \left(\frac{1}{\mu_3} \right) = 2 \sqrt{3} \quad \mbox{ and } \quad \frac{1}{\mu_3 \bar \mu_3} + \frac{1}{\mu_3^2} + \frac{1}{\bar \mu_3^2} = 2, 
$$
it follows from \Cref{lem-B} that 
\begin{multline}\label{pro-time1-p1}
\frac{1}{E}\hat u (z) \overline{\hat u (z-p)} \int_{0}^L B(z, x) \, dx 
=  \hat w (z) \overline{\hat w (z-p)} +  2 \sqrt{3} \frac{\hat w (z) \overline{\hat w (z-p)}}{(1+|z|)^{1/3}} \\[6pt] 
+  \left(2 + \frac{E_1}{E} \right) \frac{\hat w (z)  \overline{\hat w (z-p)}}{(1+|z|)^{2/3}} + \hat w (z) \overline{\hat w (z-p)} O(|z|^{-1}). 
\end{multline}

We have 
\begin{multline}\label{pro-time1-Case1-diff1}
\left| \int_{\mR} \frac{\hat w (z) \overline{\hat w (z)}}{(1+|z|)^{1/3}} \, dz  - \int_{\mR} \frac{\hat w (z) \overline{\hat w (z-p)}}{(1+|z|)^{1/3}} \, d z \right| \le  \int_{\mR} \frac{\hat w (z)}{(1 + |z|)^{1/3}} \overline{\hat w (z-p) - \hat w (z)} \, dz \\[6pt]
\le \left( \int_{\mR} \frac{|\hat w (z)|^2}{(1 + |z|)^{2/3}} \, dz  \right)^{1/2} \left(\int_{\mR} |w(t) e^{-it p} - w(t)|^2 \, dt \right)^{1/2} \\[6pt]
\mathop{\le}^{\Cref{lem-interpolation}} C T^{1/3} \| w \|_{L^2(\mR)} T \| w \|_{L^2(\mR)} = C T^{4/3} \| w\|_{L^2(\mR)}^2. 
\end{multline}
Similarly, we obtain 
\be \label{pro-time1-Case1-diff2}
\left| \int_{\mR} \frac{\hat w (z) \overline{\hat w (z)}}{(1+|z|)^{2/3}} \, dz  - \int_{\mR} \frac{\hat w (z) \overline{\hat w (z-p)}}{(1+|z|)^{2/3}} \, d z \right|  \le C T^{4/3} \| w\|_{L^2(\mR)}^2. 
\ee

Since, by \Cref{lem-interpolation}, 
$$
\int_{\mR} \frac{|\hat w (z)|^2 }{(1+|z|)^{1/3}} \, dz \le C T^{1/3} \| w\|_{L^2(\mR)}^2, 
$$
and, by H\"older's inequality,  
$$
\left| \int_{\mR} \hat w (z) \overline{\hat w (z-p)} O(|z|^{-1}) \, dz \right| \le C \int_{\mR} \frac{|\hat w(z)|^2}{1 + |z|} \, dz \mathop{\le}^{\eqref{pro-time1-01}} C \int_{|z| \le \nu^3 } \frac{|\hat w(z)|^2}{1 + |z|} \, dz \mathop{\le}^{\eqref{pro-time1-beta-nu}} C T^{3/2} \int_{\mR} |\hat w(z)|^2, 
$$
it follows from \eqref{pro-time1-p1}, \eqref{pro-time1-Case1-diff1}, and \eqref{pro-time1-Case1-diff2} that 
\begin{multline} \label{pro-time1-p2}
I : = \int_{\mR}  \frac{1}{E}\hat u (z) \overline{\hat u (z-p)} \int_{0}^L B(z, x) \, dx   \, dz \\[6pt]
= \int_{\mR} |w(t)|^2 \big(e^{- i t p} + O_+(T^{1/3}) +  O(T^{4/3})  \big) + 
\left(2 + \frac{E_1}{E} \right) \int_{\mR} \frac{|\hat w(z)|^2}{(1 + |z|^2)^{1/3}} \, dz,   
\end{multline}
where $O_+(T^\alpha)$ denotes a non-negative quantity of the order $O(T^\alpha)$ for $\alpha > 0$.

We have 
\begin{multline}
\int_{0}^T |w(t) - w(T -t)|^2 \, dt = C \int_{\mR} |\hat w(z) - e^{- i T z} \hat w(-z)|^2 \, dz \\[6pt]
= 3\int_{\mR} |\hat v_1(z) H'(z) - \hat v_1(-z) H'(-z)|^2 \, dz \\[6pt]
\mathop{\le}^{\Cref{th:lambda_asym}} C \int_{\mR} \frac{|\hat v(z)|^2  |H'(z)|^2}{(1 + |z|)^{2/3}} \, dz  
= C \int_{\mR} \frac{|\hat w(z)|^2}{(1 + |z|)^{2/3}} \, d z.   
\end{multline}
Using \Cref{lem-interpolation}, we derive that  
\be  \label{pro-time1-p3}
\int_{\mR} |w(t) - w(T -t)|^2 \, dt \le C T^{2/3} \| w\|_{L^2(\mR)}^2.
\ee

This implies 
\be \label{pro-time1-Case1-oo1}
\Re I =  \big(1 + O(T^{1/3}) \big)\int_{\mR} |w (t)|^2 \, dt , 
\ee
and, since $\Im (F/E) < 0$ by \Cref{lem-FE}, 
\be \label{pro-time1-Case1-oo2}
\Im I \mathop{\le}^{\eqref{pro-time1-01}, \eqref{pro-time1-p2}, \eqref{pro-time1-p3}} \left( \frac{ \Im(E_1/E)}{\nu^2}  - \frac{p T}{2} +  O(T^{4/3}) \right) \int_{\mR} |w (t)|^2 \, dt  . 
\ee
We derive from \eqref{pro-time1-beta-nu} 
and \Cref{lem-FE} that 
\be
\Im I \le \left( - \frac{0.1914}{6} p L  T - \frac{p T}{2} + O(T^{4/3})\right) \int_{\mR} |w (t)|^2 \, dt. 
\ee

Since, in the case $e^{i \eta_1 L } \neq 1$, 
$$
\frac{0.1914}{6} L  > 1/2 \mbox{ except for the case where } L = L(2, 1) \mbox{ or } L = L (3, 1), 
$$
it follows that, for some $\eps_0 > 0$,  
\be \label{pro-time1-Case1-oo4}
\Im I \le - (p + 2 \eps_0) T \int_{\mR} |w (t)|^2 \, dt 
\ee
when $T$ is sufficiently small.

Set 
$$
\Psi_1 = \Re \left(\frac{1}{\bar E} \Psi(T_0, x) \right) \quad \mbox{ and } \quad \Psi_2 = \Im \left(\frac{1}{\bar E} \Psi(T_0, x) \right) \mbox{ for } x \in [0, L]. 
$$

We claim that there exist $T_0 > 0$ small and $\eps_0 > 0$ small such that for all $\theta \in A$,  where 
\be
A: = \Big\{ \Psi_1 \cos (p (1+ \eps_0) t) + \Psi_2 \cos (p (1+ \eps_0) t) ; t \in [\eps_0/2, T_0]  \Big\},  
\ee
one can find a control $u_1 \in L^2([0, T_0])$ such that \footnote{Here and in what follows $\mbox{Proj}_{\cM_N} $ denotes the projection to $\cM_N$ with respect to the $L^2(0, L)$-scalar product.}
$$
y_1(T_0, \cdot) = 0 \quad \mbox{ and } \quad \mbox{Proj}_{\cM_N} y_2(T_0, \cdot) = \theta,  
$$
where $y_1$ and $y_2$ are the solutions of \eqref{eq:first_order} and \eqref{eq:second_order} with $u_2 = 0$, and zero initial datum. 

Indeed, let $T \ll T_0$,  by \eqref{projection}, \eqref{pro-time1-Case1-oo1} and \eqref{pro-time1-Case1-oo4},  there exists a control $u_{1} \in L^2(0, T_0)$ with $\supp u_1 \subset [0, T]$  such that the corresponding solution $y_1$ and $y_2$ with $u_2 = 0$ and zero initial datum satisfy  
$$
\int_0^L y_2(T_0, x) \Psi(T_0, x) = 1 + O(T^{1/3}).
$$
One then obtain, by \Cref{lem-orthogonal} that 
\be \label{pro-time1-claim-p1}
\mbox{Proj}_{\cM_N} y_2(T_0, \cdot) y_2(T_0, x) = \Psi_1(x) + O(T^{1/3}). 
\ee

On the other hand, by \eqref{projection}, \eqref{pro-time1-Case1-oo1} and \eqref{pro-time1-Case1-oo4},  there exists a control $u_{1} \in L^2([0 , T_0])$  such that the corresponding solution $y_1$ and $y_2$ with $u_2 = 0$ and zero initial datum satisfy, for some $\eps > \eps_0$,   
$y_1(T_0, \cdot) = 0$, and
$$
\int_0^L y_2(T_0, x) \Psi(T_0, x) \, dx  = e^{ - i p (1 + \eps) T_0}. 
$$
This implies 
$$
\int_0^L y_2(T_0, x) \Psi(T_0, x) e^{ - i p (1 + \eps) T_0} \, dx  = 1. 
$$
One then obtain, by \Cref{lem-orthogonal} that 
\be \label{pro-time1-claim-p2}
\mbox{Proj}_{\cM_N} y_2(T_0, \cdot) y_2(T_0, x) = \Psi_1(x) \cos \big( (1 + \eps) T_0 \big) + \Psi_2 \sin \big( (1 + \eps) T_0 \big). 
\ee

The claim now follows from \eqref{pro-time1-claim-p1} and  \eqref{pro-time1-claim-p2}.  

\medskip 
The conclusion now follows the claim using the rotation arguments as in \cite{Cerpa07}.

\medskip 
\noindent $\bullet$ Case 2: $e^{i \eta_1 L } = 1$. Fix $\gamma \neq 0$ such that 
\be \label{pro-time1-H2}
H'''(z + i \gamma) \neq 0 \mbox{ for } z \in \mR. 
\ee
Denote $H_\gamma (z) = H(z + i \gamma)$ and  set 
\be \label{pro-time1-def-uw2}
\hat u (z) =  \hat v(z) H(z) \quad \mbox{ and } \quad \hat w (z) = \frac{3^3}{\mu_3^3} z \hat v(z) H_\gamma'''(z). 
\ee
One can check that $u \in L^2(\mR)$ with $\supp u \subset [0, T]$; moreover, concerning the linearized system, $u$ is a control which steers $0$ from time $0$ to $0$ at time $T$. This follows from the fact that  $
e^{i \beta \cdot}\hu$ and $ e^{i \beta \cdot}  \hu G/H$ are entire functions and satisfy the assumptions of Paley-Wiener's theorem concerning the support in $[-T/2, T/2]$ (see also \Cref{lem-form-sol} and \Cref{pro-Gen}).

We have, by \Cref{th:lambda_asym} (see also \eqref{pro-time1-recall}), for real $z$ with large $|z|$, 
$$
\hat w(z) = \frac{\hat u (z)}{|z|} \left(1 - \frac{6}{\mu_3 |z|^{1/3}} +  \frac{31}{\mu_3^2 |z|^{2/3}} + O(|z|^{-1}))  \right).  
$$
This implies, by \eqref{pro-time1-H2}, for $z \in \mR$, 
$$
\frac{\hat u (z)}{|z|} =\hat w(z) \Big(1 + \frac{6}{\mu_3 |z|^{1/3}}  + \frac{5}{\mu_3^2 |z|^{2/3}} + 
O (|z|^{-1}) \Big). 
$$
We thus derive that 
\begin{multline} \label{pro-time1-Case2-p1}
\hat u (z) \overline{\hat u (z-p)} \left( \frac{1}{|z|^{2}} + \frac{F_1}{F} \frac{1}{|z|^{8/3}} + O(|z|^{-9/3}) \right)
= \hat w (z) \overline{\hat w (z-p)} + 2 \Re\left( \frac{6}{\mu_3} \right) \frac{\hat w (z) \overline{\hat w (z-p)}}{|z|^{1/3}} \\[6pt]
+  \left( \frac{F_1}{F} + 5 \left( \frac{1}{\mu_3^2} + \frac{1}{\mu_3 \tmu_3}  + \frac{1}{\tmu_3^2}  \right) \right)  \frac{\hat w (z)  \overline{\hat w (z-p)}}{|z|^{2/3}} + \hat w (z) \overline{\hat w (z-p)} O(|z|^{-1}). 
\end{multline}
Since 
$$
2 \Re \left( \frac{6}{\mu_3 |z|^{1/3}} \right) = 6 \sqrt{3} \quad \mbox{ and } \quad \frac{1}{\mu_3^2} + \frac{1}{\mu_3 \tmu_3}  + \frac{1}{\tmu_3^2}  =2, 
$$
it follows from  \Cref{lem-B} that  
\begin{multline}
\frac{1}{F} \hat u (z) \overline{\hat u (z-p)} \int_0^L B(z, x) \, dx  \, dz 
= \hat w (z) \overline{\hat w (z-p)} \\[6pt] + 6 \sqrt{3} \frac{\hat w (z) \overline{\hat w (z-p)} }{|z|^{1/3}} + 
\left( \frac{F_1}{F} + 10 \right)  \frac{\hat w (z)  \overline{\hat w (z-p)}}{|z|^{2/3}} + \hat w (z) \overline{\hat w (z-p)} O(|z|^{-1})
\end{multline}

Similar to \eqref{pro-time1-p2}, we obtain 
\begin{multline} \label{pro-time1-Case2-p2}
J: = \int_{\mR} \frac{1}{F} \hat u (z) \overline{\hat u (z-p)} \int_0^L B(z, x) \, dx  \, dz \\[6pt]
=   \int_{\mR} |w(t)|^2 \big(e^{- i t p} +  O_+(T^{1/3}) + O(T^{4/3}) \big) + \left(\frac{F_1}{F } +  10 \right)\int_{\mR} \frac{|\hat w(z)|^2}{(1 + |z|)^{2/3}}\, dz.
\end{multline}
We have, by \Cref{lem-B,lem-G}, 
$$
\frac{F_1}{F} = - \frac{2}{3} - \frac{ipL}{3}. 
$$

Similar to \eqref{pro-time1-Case1-oo1} and \eqref{pro-time1-Case1-oo2}, we obtain 
\be \label{pro-time1-Case2-oo1}
\Re J = (1 + O(T)) \int_{\mR} |w (t)|^2 \, dt, 
\ee
and
\be \label{pro-time1-Case2-oo4}
\Im J \le \left( 0.1914 \Im(F_1/ F) T  - \frac{pT}{2}  \right) \int_{\mR} |w (t)|^2 \, dt  \le \left(-  \frac{0.1914 L}{3} - \frac{1}{2} \right)p T \int_{\mR} |w (t)|^2 \, dt . 
\ee
We have 
$$
\frac{0.1914 L}{3} \ge \frac{0.1914 L(2, 2)}{3}  \ge 0.801 > 1/2,   
$$

Set 
$$
\Psi_1 = \Re \left(\frac{1}{\bar E} \Psi(T_0, x) \right) \quad \mbox{ and } \quad \Psi_2 = \Im \left(\frac{1}{\bar E} \Psi(T_0, x) \right) \mbox{ for } x \in [0, L]. 
$$
We claim that there exist $T_0 > 0$ small and $\eps_0 > 0$ small such that for all $\theta \in A$,  where 
\be
A: = \Big\{ \Psi_1 \cos (p (1+ \eps_0) t) + \Psi_2 \cos (p (1+ \eps_0) t) ; t \in [\eps_0/2, T_0]  \Big\},  
\ee
one can find a control $u_1 \in L^2([0, T_0])$ such that 
$$
y_1(T_0, \cdot) = 0 \quad \mbox{ and } \quad \mbox{Proj}_{\cM_N} y_2(T_0, \cdot) = \theta,  
$$
where $y_1$ and $y_2$ are the solutions of \eqref{eq:first_order} and \eqref{eq:second_order} with $u_2 = 0$, and zero initial datum. 
The proof of this claim is similar to the corresponding one in previous case using  \eqref{pro-time1-Case2-oo1} and \eqref{pro-time1-Case2-oo4}, and \Cref{lem-orthogonal}. The proof is omitted. 

\medskip 
The conclusion now follows the claim using the rotation arguments as in \cite{Cerpa07}. 
\end{proof}

\subsection{Proof of \Cref{thm-time}} \Cref{thm-time} now follows from \Cref{pro-time1} and the arguments of Cerpa and Cr\'epeau \cite{CC09}. The details are omitted to avoid the repetition.    

\appendix
\section{A lemma on interpolation inequalities}

\renewcommand{\theequation}{A\arabic{equation}}
\renewcommand{\theproposition}{A\arabic{proposition}}
\renewcommand{\thelemma}{A\arabic{lemma}}
  \setcounter{equation}{0}  
  \setcounter{lemma}{0}  

\begin{lemma} \label{lem-interpolationA} Let $T > 0$. There exists a positive constant $C$ such that  
\be \label{lem-interpolation1-cl1}
\| u\|_{L^2(0, T)} \le C \| u\|_{[H^{2/3}(0, T)]^*}^{3/7} \| u\|_{H^{1/2}(0, T)}^{4/7} \mbox{ for } u \in H^{1/2}(0, T),  
\ee
\be \label{lem-interpolation1-cl1-2}
\| u\|_{[H^{1/3}(0, T/2)]^*} \le C \| u\|_{[H^{2/3}(0, T/2)]^*}^{5/7} \| u\|_{H^{1/2}(0, T/2)}^{2/7}  \mbox{ for } u \in H^{1/2}(0, T), 
\ee
\be \label{lem-interpolation1-cl2}
\| u\|_{L^2(0, T)} \le C \| u\|_{[H^{1}(0, T)]^*}^{7/13} \| u\|_{H^{7/6}(0, T)}^{6/13} \mbox{ for } u \in H^{7/6}(0, T),  
\ee
and 
\be \label{lem-interpolation1-cl2-2}
\| u\|_{[H^{1/3}(0, T/2)]^*} \le C \| u\|_{[H^{1}(0, T/2)]^*}^{9/13} \| u\|_{H^{7/6}(0, T/2)}^{4/13} \mbox{ for } u \in H^{7/6}(0, T).  
\ee
\end{lemma}

\begin{proof} We first note that 

\be \label{lem-interpolation1-p1}
\| v\|_{L^2(\mR)} \le C \| v \|_{H^{-2/3}(\mR)}^{3/7} \| v\|_{H^{1/2}(\mR)}^{4/7} \mbox{ for } v \in H^{1/2}(\mR). 
\ee
\be \label{lem-interpolation1-p1-2}
\| v\|_{H^{-1/3}(\mR)} \le C \| v\|_{H^{-2/3}(\mR)}^{5/7} \| v\|_{H^{1/2}(\mR)}^{2/7}  \mbox{ for } u \in H^{1/2}(0, T), 
\ee 
\be \label{lem-interpolation1-p2}
\| v\|_{L^2(\mR)} \le C \| v\|_{H^{-1}(\mR)}^{7/13} \| v\|_{H^{7/6}(\mR)}^{6/13} \mbox{ for } v \in H^{7/6}(\mR),  
\ee
and
\be \label{lem-interpolation1-p2-2}
\| v\|_{H^{-1/3}(\mR)} \le C \| v \|_{H^{-1}(\mR)}^{9/13} \| v \|_{H^{7/6}(\mR)}^{4/13} \mbox{ for } u \in H^{7/6}(\mR).  
\ee
These inequalities follow from the fact that, by H\"older's inequality, 
\begin{multline*}
\int_{\mR} |\hat v(\xi)|^2 \, d \xi  =  \int_{\mR} \frac{|\hat v(\xi)|^{\frac{6}{7}}}{(1 + |\xi|^2)^{\frac{2}{7}}}  |\hat v(\xi)|^{\frac{8}{7}}(1 + |\xi|^2)^{\frac{2}{7}} \, d \xi \\[6pt]
\le \left( \int_{\mR} \frac{|\hat v(\xi)|^{2}}{(1 + |\xi|^2)^{\frac{2}{3}}} \, d \xi \right)^{\frac{3}{7}} \left( \int_{\mR} |\hat v(\xi)|^{2}(1 + |\xi|^2)^{\frac{1}{2}} \, d \xi \right)^{\frac{4}{7}}. 
\end{multline*}
\begin{multline*}
\int_{\mR} \frac{|\hat v(\xi)|^2}{(1 + |\xi|^2)^{1/3}} \, d \xi  =  \int_{\mR} \frac{|\hat v(\xi)|^{\frac{10}{7}}}{(1 + |\xi|^2)^{\frac{10}{21}}}  |\hat v(\xi)|^{\frac{4}{7}}(1 + |\xi|^2)^{\frac{3}{21}} \, d \xi \\[6pt]
\le \left( \int_{\mR} \frac{|\hat v(\xi)|^{2}}{(1 + |\xi|^2)^{\frac{2}{3}}} \, d \xi \right)^{\frac{5}{7}} \left( \int_{\mR} |\hat v(\xi)|^{2}(1 + |\xi|^2)^{\frac{1}{2}} \, d \xi \right)^{\frac{2}{7}}, 
\end{multline*}
\begin{multline*}
\int_{\mR} |\hat v(\xi)|^2 \, d \xi  =  \int_{\mR} \frac{|\hat v(\xi)|^{\frac{14}{13}}}{(1 + |\xi|^2)^{\frac{7}{13}}}  |\hat v(\xi)|^{\frac{15}{13}}(1 + |\xi|^2)^{\frac{7}{13}} \, d \xi \\[6pt]
\le \left( \int_{\mR} \frac{|\hat v(\xi)|^{2}}{(1 + |\xi|^2)} \, d \xi \right)^{\frac{7}{13}} \left( \int_{\mR} |\hat v(\xi)|^{2}(1 + |\xi|^2)^{\frac{7}{6}} \, d \xi \right)^{\frac{6}{13}}. 
\end{multline*}
and 
\begin{multline*}
\int_{\mR} \frac{|\hat v(\xi)|^2}{(1 + |\xi|^2)^{1/3}}  \, d \xi  =  \int_{\mR} \frac{|\hat v(\xi)|^{\frac{18}{13}}}{(1 + |\xi|^2)^{\frac{9}{13}}}  |\hat v(\xi)|^{\frac{8}{13}}(1 + |\xi|^2)^{\frac{14}{39}} \, d \xi \\[6pt]
\le \left( \int_{\mR} \frac{|\hat v(\xi)|^{2}}{(1 + |\xi|^2)} \, d \xi \right)^{\frac{9}{13}} \left( \int_{\mR} |\hat v(\xi)|^{2}(1 + |\xi|^2)^{\frac{7}{6}} \, d \xi \right)^{\frac{4}{13}}. 
\end{multline*}

Let $w$ be a function defined in $[-T, 2T]$ by   
\be
w(t) = \left\{\begin{array}{cl} u(t) &  \mbox{ in } (0, T) \\[6pt]
- 3 u(-t) + 4 u(-t/2) &  \mbox{ in } (-T, 0), \\[6pt]
-3 u(2 T - t) + 4 u \left( \frac{3 T - t}{2} \right) & \mbox{ in } (T, 2T).  
\end{array}\right. 
\ee
There is a constant $\Lambda \ge 1$ independent of $u$, such that, for $0 \le s_1 < 2$, and for $0 < s_2 < 1$,     
$$
\| u \|_{H^{s_1}(0, T)} \le  \| w \|_{H^{s_1}(-T, 2 T)} \le \Lambda \| u \|_{H^{s_1}(0, T)},   
$$
and 
$$
\| u \|_{[H^{s_2}(0, T)]^*} \le \| w \|_{[H^{s_2}(-T, 2 T)]^*} \le \Lambda \| u \|_{[H^{s_2}(0, T)]^*}.  
$$
Fix $\varphi  \in C^\infty_{c}(-T, 2 T)$ such that $\varphi = 1$ in $[0, T]$ and $\varphi = 0 $ in $[-T, - T/2] \cup [3T/2, 2 T]$ and set 
$$
v = \varphi w \mbox{ in } \mR. 
$$
Then, for some $C \ge 1$ independent of $u$,  
$$
C^{-1}\| u \|_{H^{s_1}(0, T)} \le  \| v \|_{H^{s_1}(\mR)} \le C \| u \|_{H^{s_1}(0, T)},   
$$
and 
$$
C^{-1} \| u \|_{[H^{s_2}(0, T)]^*} \le \| v \|_{H^{-s_2}(\mR)} \le C \| u \|_{[H^{s_2}(0, T)]^*}.  
$$
Assertions \eqref{lem-interpolation1-cl1}, \eqref{lem-interpolation1-cl1-2},   \eqref{lem-interpolation1-cl2}, and  \eqref{lem-interpolation1-cl2-2} now follow from \eqref{lem-interpolation1-p1}, \eqref{lem-interpolation1-p1-2},  \eqref{lem-interpolation1-p2}, and \eqref{lem-interpolation1-p2-2}, respectively.  
\end{proof}

\end{document}